\documentclass[a4paper,12pt]{amsart}

\allowdisplaybreaks
\usepackage[T1]{fontenc}
\usepackage[utf8]{inputenc}

\usepackage{amsmath}
\usepackage{amsthm}
\usepackage{amssymb}
\usepackage{amstext}
\usepackage{mathtools}

\usepackage{mathrsfs}
\usepackage{bm}
\usepackage{shuffle}
\usepackage[all]{xy}
\usepackage{ytableau}
\ytableausetup{boxsize=normal,aligntableaux=center}

\usepackage{xcolor} 
\usepackage{graphicx}
\usepackage{tikz}

\usepackage{ascmac}
\usepackage{comment}
\usepackage{here}
\usepackage{time}
\usepackage{fullpage} 

\usepackage{hyperref}
\hypersetup{
  colorlinks=true,
  linkcolor=blue,
  citecolor=blue,
  filecolor=blue,
  urlcolor=blue
}

\usepackage[capitalize,nameinlink,noabbrev,nosort]{cleveref}

\DeclareFontEncoding{OT2}{}{}
\DeclareFontSubstitution{OT2}{cmr}{m}{l}
\DeclareFontFamily{OT2}{cmr}{\hyphenchar\font45 }
\DeclareFontShape{OT2}{cmr}{m}{l}{%
  <5><6><7><8><9>gen*wncyr%
  <10><10.95><12><14.4><17.28><20.74><24.88>wncyr10}{}
\DeclareMathAlphabet{\mathcyr}{OT2}{cmr}{m}{l}
\DeclareMathAlphabet{\mathcyb}{OT2}{cmr}{b}{l}
\SetMathAlphabet{\mathcyr}{bold}{OT2}{cmr}{b}{l}

\makeatletter
\@namedef{subjclassname@2020}{%
  \textup{2020} Mathematics Subject Classification}
\makeatother

\newtheorem{theoremcounter}{Theorem Counter}[section]

\theoremstyle{definition}
\newtheorem{definition}[theoremcounter]{Definition}
\newtheorem{remark}[theoremcounter]{Remark}
\newtheorem{example}[theoremcounter]{Example}

\theoremstyle{plain}
\newtheorem{lemma}[theoremcounter]{Lemma}
\newtheorem{proposition}[theoremcounter]{Proposition}
\newtheorem{corollary}[theoremcounter]{Corollary}
\newtheorem{conjecture}[theoremcounter]{Conjecture}
\newtheorem{theorem}[theoremcounter]{Theorem}

\theoremstyle{definition}

\theoremstyle{plain}

\numberwithin{equation}{section}

\newenvironment{nouppercase}{%
  \renewcommand{\uppercasenonmath}[1]{}}{}
\begin{document}

\title{Anti-ribbon shaped Young diagram decomposition}

\author[H. Murahara]{Hideki Murahara}
\address{The University of Kitakyushu,
4-2-1 Kitagata, Kokuraminami-ku,
Kitakyushu, Fukuoka 802-8577,
Japan}
\email{hmurahara@mathformula.page}

\author{Maki Nakasuji}
\address{Department of Information and Communication Science, Faculty of Science, Sophia
University, 7-1 Kioi-cho, Chiyoda-ku, Tokyo, 102-8554, Japan}
\email{nakasuji@sophia.ac.jp}

\author[W. Takeda]{Wataru Takeda}
\address{Department of Mathematics, Toho University,
2-2-1 Miyama, Funabashi, Chiba 274-8510, Japan.}
\email{wataru.takeda@sci.toho-u.ac.jp}

\subjclass[2020]{Primary 05E05; Secondary 11M32}
\keywords{Sliding procedure, Anti-ribbon shape,  Schur function, Schur multiple zeta function}

\begin{abstract}
We introduce a sliding procedure on diagrams which splits the standard Young tableaux of a corner-crossed diagram, in particular of a skew shape, along a finite binary tree whose leaves are anti-ribbon diagrams.
Combined with standardization, this yields a decomposition of semi-standard Young tableaux into tableaux of anti-ribbon shape, indexed by standard tableaux and governed by their descent sets.
As applications, we obtain expansions of Schur functions and Schur multiple zeta functions in terms of their anti-ribbon type analogues; at the level of Schur functions, this is a tableau-level form of Gessel's fundamental quasisymmetric expansion.
\end{abstract}

\begin{nouppercase}
\maketitle
\end{nouppercase}

\section{Introduction}
Schur functions, which play a central role in representation theory and combinatorics, are well known to admit a combinatorial description via their tableau formula. 
This expression is given as a generating series indexed by semi-standard Young tableaux whose shapes are determined by partitions, that is, weakly decreasing sequences of positive integers.

On the other hand, as an analogue of the tableau expression for Schur functions, Nakasuji, Phuksuwan, and Yamasaki~\cite{NakasujiPhuksuwanYamasaki2018} introduced Schur multiple zeta functions as a combinatorial generalization of multiple zeta functions.
These functions are defined as multivariate series indexed by semi-standard Young tableaux, in which the tableau entries serve as running indices. 
By construction, Schur multiple zeta functions can be expressed as linear combinations of classical multiple zeta functions. 
Consequently, the problem of describing a basis for the space spanned by Schur multiple zeta values is equivalent to that for multiple zeta values themselves.

In this study, we introduce a procedure, termed \emph{sliding}, which acts on the boxes of a diagram.
An anti-ribbon shape is a connected diagram containing no $2\times 2$ block and running from the upper-left to the lower-right; unlike an ordinary Young diagram, it is naturally parametrized by a composition rather than by a partition.
We prove that the sliding procedure, applied to a corner-crossed diagram, that is, a diagram whose northwest corners are pairwise crossing, terminates after finitely many steps, and that every terminal diagram is a single anti-ribbon diagram (Corollaries~\ref{cor:termination} and~\ref{cor:leaf-antiribbon}); this class contains all Young diagrams, all skew Young diagrams and all edge-connected diagrams.
The resulting finite binary tree has exactly $|\mathrm{SYT}(D)|$ leaves, and the slides along each root-to-leaf path arrange the boxes of $D$ in the order prescribed by a unique standard tableau (Theorem~\ref{thm:onetoone}).
Combined with the standardization of semi-standard tableaux, this partitions the semi-standard tableaux of a fixed corner-crossed diagram, in particular of a fixed skew shape, into blocks indexed by the standard tableaux, each block being in bijection with the semi-standard tableaux of an anti-ribbon shape determined by the descent set of the corresponding standard tableau (Theorem~\ref{thm:standardization}).

A Schur function and a Schur multiple zeta function are each defined as a sum over the semi-standard tableaux $M=(m_{ij})$ of a fixed shape. The summand is the monomial $\prod x_{m_{ij}}$ in the first case, and the product of powers $\prod m_{ij}^{-s_{ij}}$ in the second, where the weight $\pmb s=(s_{ij})$ is a tableau of complex variables on the same shape.
A slide only relocates boxes, without duplicating or discarding them. Hence, it leaves the monomial $\prod x_{m_{ij}}$ unchanged, and it carries the weight $s_{ij}$ along with the box to which it is attached, so that the summand of the Schur multiple zeta function is unchanged as well.
As a consequence, the above bijection turns into a linear expansion of the Schur function and of the Schur multiple zeta function as a sum of the corresponding anti-ribbon type functions attached to the leaves.
We note that the determinant and Pfaffian expressions of \cite{NakasujiPhuksuwanYamasaki2018,NakasujiTakeda2022,MatsumotoNakasuji2024,NakasujiTakeda2025}, see also \cite{TakedaYamasaki2025} for algebraic relations derived from them,
express a Schur multiple zeta function through functions of hook or of one-row type, whereas the expansion obtained here is a sum, with no signs and no cancellation, over anti-ribbon shapes.
At the level of Schur functions, the resulting identity is Gessel's expansion of a skew Schur function into fundamental quasisymmetric functions~\cite{Gessel1984}.
Reading the variables of a Schur multiple zeta function along the standard tableaux of its shape also underlies the expression of $\zeta_\lambda$ as a sum of multiple zeta functions in \cite[Section~2.2]{NakasujiPhuksuwanYamasaki2018}; here the contributions attached to one standard tableau are collected into a single Schur multiple zeta function of anti-ribbon type, and the order in which the boxes are read is produced by the sliding procedure.

The paper is organized as follows.
In Section~\ref{sec:prelim}, we fix the notation concerning partitions, skew Young diagrams and tableaux, and collect the dictionary relating anti-ribbon diagrams, compositions and $rb$-words which is used throughout.
In Section~\ref{sec:sliding}, we define the sliding procedure on diagrams, establish its basic properties, and prove the bijection described above.
In Section~\ref{subsec:anti-ribbon-descents}, we describe the anti-ribbon shapes in terms of descent words of standard tableaux, prove the standardization theorem, and compare the result with the terminal shapes of the sliding tree.
In Section~\ref{sec:applications}, we apply these results to Schur functions and Schur multiple zeta functions.

\section{Preliminaries}
\label{sec:prelim}
\subsection{Diagrams and tableaux}
Let $\lambda=(\lambda_1,\ldots,\lambda_k)$ be a partition, that is, a weakly decreasing sequence of positive integers.
We associate with $\lambda$ the Young diagram
\[
D(\lambda)
=\{(i,j)\in \mathbb{Z}_{>0}^2
\mid
1\leq i\leq k,\ 1\leq j\leq \lambda_i\}.
\]
Thus the $i$-th row of $D(\lambda)$ consists of $\lambda_i$ boxes.
Let $\lambda=(\lambda_1, \ldots, \lambda_k)$ and $\mu=(\mu_1, \ldots, \mu_k)$ be partitions such that $\mu\subset \lambda$, which means $\mu_i \leq \lambda_i$ for all $i$.
The skew Young diagram of shape $\lambda/\mu$ is defined by
$D(\lambda/\mu)=\{(i, j)\in {\mathbb Z}^2 \mid 1\leq i \leq k, \mu_i<j\leq \lambda_i \}$,
which is obtained by removing the Young diagram 
$D(\mu)$ from $D(\lambda)$; its $i$-th row contains $\lambda_i-\mu_i$ boxes.
We say that $(i, j)\in D(\lambda/\mu)$ is a corner of $\lambda/\mu$ if $(i+1,j)\not\in D(\lambda/\mu)$ and $(i,j+1)\not\in D(\lambda/\mu)$, and we denote by $C(\lambda/\mu)\subset D(\lambda/\mu)$ the set of all corners.
For the standard theory of partitions, skew Young diagrams, Schur functions and ribbon shapes we refer to Macdonald~\cite{Macdonald1995}.

More generally, by a \emph{diagram} we mean a finite subset $D\subset\mathbb{Z}_{>0}^2$, the pair $(i,j)$ being regarded as the box in the $i$-th row and the $j$-th column, where row indices increase from top to bottom and column indices increase from left to right; a skew Young diagram is a diagram in this sense.
A tableau on $D$ over a set $X$ is a filling $(t_{ij})$ of $D$ with entries $t_{ij}\in X$; we write $T(D,X)$ for the set of all of them and $T(D)=T(D,\mathbb{N})$, where $\mathbb{N}=\mathbb{Z}_{\ge1}$.
We abbreviate
\[
T(\lambda/\mu,X)=T\bigl(D(\lambda/\mu),X\bigr). 
\]

The diagrams produced by the sliding procedure of Section~\ref{sec:sliding} may have gaps inside a row or a column.
The order relation used throughout the paper is therefore defined by means of the following auxiliary cells.
\begin{definition}[Implicit boxes and the row--column order]
\label{def:implicit}
Let $D$ be a diagram and let $\widetilde{D}$ be the smallest set of cells containing $D$ and closed under the following two operations: if two cells of $\widetilde{D}$ lie in the same row, every cell between them in that row belongs to $\widetilde{D}$; if two cells of $\widetilde{D}$ lie in the same column, every cell between them in that column belongs to $\widetilde{D}$. The cells of $\widetilde{D}\setminus D$ are called \emph{implicit boxes}; they carry no entry.

The \emph{row--column order} $\preceq_D$ on $D$ is the restriction to $D$ of the row--column order of $\widetilde{D}$: for $u,v\in D$ we set $u\preceq_D v$ if there is a chain from $u$ to $v$ in $\widetilde{D}$ each step of which goes to the right within a row or downward within a column.
\end{definition}
If $D$ has no implicit box, in particular if $D$ is a Young or a skew Young diagram, then $\preceq_D$ is simply the transitive closure of the relations
\[
(i,j)\preceq_D(i,j')\quad (j\le j'),
\qquad
(i,j)\preceq_D(i',j)\quad (i\le i'),
\]
taken over the boxes belonging to $D$. In general the two prescriptions differ: for
\[
D=\{(1,1),(1,2),(2,1),(2,2),(3,1),(3,3)\}
\]
the cell $(3,2)$ is implicit, and one has $(1,2)\preceq_D(2,2)\preceq_D(3,3)$, whereas $(1,2)$ and $(3,3)$ would be incomparable if the implicit box were disregarded.

A filling $T\in T(D)$ is a \emph{standard Young tableau} of $D$, if it is a bijection from $D$ onto $\{1,\ldots,|D|\}$ such that
\[
x\preceq_D y\ \Longrightarrow\ T(x)\leq T(y);
\]
we denote by $\mathrm{SYT}(D)$ the set of all standard Young tableaux on $D$, and abbreviate $\mathrm{SYT}(\lambda/\mu)=\mathrm{SYT}(D(\lambda/\mu))$.
If $D$ has no implicit box, this condition says exactly that the entries increase along each row from left to right and along each column from top to bottom, so that $\mathrm{SYT}(\lambda/\mu)$ is the usual set of standard tableaux of shape $\lambda/\mu$.
In general the relations imposed through implicit boxes are strictly stronger, and it is these that are needed in Section~\ref{sec:sliding}.
More generally, for a diagram $D$ we denote by $\mathrm{SSYT}(D)$ the set of those $M\in T(D)$ such that, for $u=(i,j)$ and $v=(i',j')$ in $D$,
\[
u\preceq_D v\ \Longrightarrow\ M(u)\le M(v),
\qquad
u\preceq_D v,\ i<i'\ \Longrightarrow\ M(u)<M(v).
\]
We abbreviate
\[
\mathrm{SSYT}(\lambda/\mu)=\mathrm{SSYT}\bigl(D(\lambda/\mu)\bigr).
\]
If $D$ has no implicit box, in particular if $D$ is a skew Young or an anti-ribbon diagram, this condition says exactly that the entries are weakly increasing along each row and strictly increasing along each column: a chain realizing $u\preceq_D v$ then runs through boxes of $D$, the entries weakly increase along its horizontal steps and strictly along its vertical ones, and a vertical step occurs when $i<i'$.

\subsection{Anti-ribbon diagrams, words and compositions}
\label{subsec:antiribbon}
A ribbon shape is a connected skew diagram with no $2\times 2$ block, arranged along a path from the upper-right to the lower-left. 
We are concerned with the reverse configuration: a connected diagram, with no $2\times2$ block arranged along a path from the upper-left to the lower-right is termed an \emph{anti-ribbon} diagram, and a tableau whose shape is such a diagram is termed an anti-ribbon tableau. 
Such shapes are described by compositions: for a composition $\alpha = (\alpha_1, \ldots, \alpha_k)$, the anti-ribbon diagram $D_{\mathrm{ar}}(\alpha)$ is the diagram consisting of $k$ rows, where the $i$-th row has $\alpha_i$ boxes and, for each $i = 1, \ldots, k-1$, the leftmost box of the $(i+1)$-st row is placed immediately below the rightmost box of the $i$-th row. 
Equivalently, $D_{\mathrm{ar}}(\alpha)$ is obtained by starting from the upper-left box and moving repeatedly either one step to the right or one step downward; we refer to this description as the \emph{southeast convention} and use it throughout.

It is convenient to encode anti-ribbon shapes by words as well.
For $N\ge1$ put
\[
\mathcal{W}_N
=
\{\,s\varepsilon_1\varepsilon_2\cdots \varepsilon_{N-1}
\mid
\varepsilon_k\in\{r,b\}\,\},
\]
and read the letter $r$ as a step to the right and the letter $b$ as a step downward.

\begin{definition}[Anti-ribbon diagram attached to a word]
\label{def:Rw}
Let
\[
w=s\varepsilon_1\varepsilon_2\cdots \varepsilon_{N-1}\in\mathcal{W}_N .
\]
Set $p_1=(1,1)$ and, once $p_k=(i,j)$ has been defined,
\[
p_{k+1}
=
\begin{cases}
(i,j+1),& \varepsilon_k=r,\\
(i+1,j),& \varepsilon_k=b,
\end{cases}
\qquad
R(w)=\{p_1,p_2,\ldots,p_N\}.
\]
Thus $R(w)$ is an anti-ribbon diagram, and $p_1,\ldots,p_N$ are its boxes in path order.
\end{definition}

These three encodings, by a composition, by a word, and by a subset of $\{1,\ldots,N-1\}$, correspond to one another as follows.

\begin{definition}
\label{def:dictionary}
For $w=s\varepsilon_1\cdots\varepsilon_{N-1}\in\mathcal{W}_N$ put
\[
\operatorname{Des}(w)=\{\,k\in\{1,\ldots,N-1\}\mid \varepsilon_k=b\,\},
\]
and, writing $\operatorname{Des}(w)=\{d_1<\cdots<d_\ell\}$,
\[
\alpha(w)
=
(d_1,\,d_2-d_1,\,\ldots,\,d_\ell-d_{\ell-1},\,N-d_\ell),
\]
with $\alpha(w)=(N)$ when $\operatorname{Des}(w)=\emptyset$.
Conversely, for a composition $\alpha=(\alpha_1,\ldots,\alpha_m)$ of $N$ put
\[
S(\alpha)=\{\alpha_1,\ \alpha_1+\alpha_2,\ \ldots,\ \alpha_1+\cdots+\alpha_{m-1}\}\subset\{1,\ldots,N-1\}.
\]
\end{definition}

\begin{lemma}
\label{lem:dictionary}

For $w\in \mathcal{W}_N$, the maps $w\mapsto\operatorname{Des}(w)$, $w\mapsto\alpha(w)$, and $\alpha\mapsto S(\alpha)$ are bijections from $\mathcal{W}_N$ onto the subsets of $\{1,\ldots,N-1\}$, from $\mathcal{W}_N$ onto the compositions of $N$, and from the compositions of $N$ onto the subsets of $\{1,\ldots,N-1\}$, respectively; moreover
\[
S(\alpha(w))=\operatorname{Des}(w),
\qquad
R(w)=D_{\mathrm{ar}}(\alpha(w))
\]
(the latter up to translation).
\end{lemma}

\begin{proof}
A word in $\mathcal{W}_N$ is determined by the positions of its letters $b$, that is, by $\operatorname{Des}(w)\subseteq\{1,\ldots,N-1\}$.
Then, a subset $\{d_1<\cdots<d_\ell\}$ of $\{1,\ldots,N-1\}$ corresponds to the composition of $N$ formed by the successive gaps, which is exactly $\alpha(w)$.
This is inverse to $\alpha\mapsto S(\alpha)$.
In $R(w)$, a maximal run of letters $r$ produces a maximal horizontal run of boxes, and each letter $b$ starts a new row directly below the last box of the previous one.
Thus, the row lengths are the successive gaps of $\operatorname{Des}(w)$, which is the description of $D_{\mathrm{ar}}(\alpha(w))$ given above.
\end{proof}

We write $T_{\mathrm{ar}}(\alpha,X)=T(D_{\mathrm{ar}}(\alpha),X)$ and $\mathrm{SSYT}_{\mathrm{ar}}(\alpha)=\mathrm{SSYT}(D_{\mathrm{ar}}(\alpha))$.
Semi-standardness on an anti-ribbon shape is a condition on consecutive boxes only.

\begin{lemma}
\label{lem:ar-ssyt}
Let $\alpha$ be a composition of $N$ and let $p_1,\ldots,p_N$ be the boxes of $D_{\mathrm{ar}}(\alpha)$ in path order.
A filling of $D_{\mathrm{ar}}(\alpha)$ whose entry at $p_k$ is $i_k$ is semi-standard if and only if
\[
i_1\le i_2\le\cdots\le i_N
\qquad\text{and}\qquad
i_t<i_{t+1}\quad (t\in S(\alpha)).
\]
\end{lemma}

\begin{proof}
Consecutive boxes $p_t,p_{t+1}$ lie in the same row when $t\notin S(\alpha)$ and in the same column when $t\in S(\alpha)$, so the displayed conditions are exactly the monotonicity conditions for consecutive boxes.
Conversely, the boxes of $D_{\mathrm{ar}}(\alpha)$ lying in a fixed row, respectively in a fixed column, form a set of consecutive terms of the sequence $p_1,\ldots,p_N$.
Thus, the general monotonicity conditions follow from those for consecutive boxes by transitivity, and we obtain the lemma.
\end{proof}

\section{The sliding procedure}
\label{sec:sliding}
We give a diagram-based formulation of the sliding procedure. It is expressed in terms of outer pre-ordered boxes, and it naturally produces a rooted binary tree.

\begin{remark}
\label{rem:jdt}
The word slide is also used for the elementary moves of the jeu de taquin of Sch\"utzenberger~\cite{Schutzenberger1977}.
A jeu de taquin slide moves an entry into an adjacent empty cell, and its iteration rectifies a skew tableau to a single tableau of straight shape. The slides of Definition~\ref{def:slides} instead move a box to destroy a crossing in the sense of Definition~\ref{def:crossing}: two boxes are moved in the two possible ways at each step, producing a binary tree of anti-ribbon shapes rather than a single normal form.
\end{remark}

\subsection{Slides and the sliding tree}
\begin{definition}[Crossings and pre-ordered boxes]
\label{def:crossing}
Let $D$ be a diagram. An unordered pair $\{u,v\}$ of distinct boxes of $D$, say $u=(i,j)$ and $v=(i',j')$, is termed a \emph{crossing} of $D$ if $v$ lies in the first or the third open quadrant of the $xy$-plane centered at $u$, that is, if
\[
(i'-i)(j-j')>0.
\]
We write $\operatorname{Cr}(D)$ for the set of crossings of $D$. Crossings are defined using the boxes of $D$ only; the implicit boxes of Definition~\ref{def:implicit} play no role here.

A box is \emph{pre-ordered} if it belongs to some crossing, and we write $\operatorname{Pre}(D)$ for the set of pre-ordered boxes; thus $\operatorname{Pre}(D)=\emptyset$ if and only if $\operatorname{Cr}(D)=\emptyset$.
An \emph{outer pre-ordered box} is a pre-ordered box which is maximal for the row--column order $\preceq_D$ of Definition~\ref{def:implicit} among the pre-ordered boxes, that is, a box $b\in\operatorname{Pre}(D)$ such that
\[
v\in\operatorname{Pre}(D),\quad b\preceq_D v
\quad\Longrightarrow\quad
v=b .
\]
We write $\operatorname{Out}(D)$ for the set of outer pre-ordered boxes.
Since two boxes of $D$ lying in one row or in one column are comparable for $\preceq_D$, outer pre-ordered boxes have pairwise distinct row indices and pairwise distinct column indices.
We may therefore always list the outer pre-ordered boxes as $b^{(1)},\ldots,b^{(m)}$ in increasing order of their rows.
\end{definition}

\begin{definition}[Admissible pair]
\label{def:admissible}
A pair $(b_1,b_2)$ of boxes of $D$ is \emph{admissible} if $b_1=(p,q)\in\operatorname{Out}(D)$, $b_2=(r,s)\in\operatorname{Out}(D)$ and
\[
p<r,\qquad q>s,
\]
that is, $b_2$ lies strictly to the lower-left of $b_1$.
With the outer pre-ordered boxes listed as $b^{(1)},\ldots,b^{(m)}$ in increasing order of their rows, the pair $(b_1,b_2)$ is called \emph{consecutive} if $b_1=b^{(t)}$ and $b_2=b^{(t+1)}$ for some $t$.
\end{definition}
We now define two sliding operations. Throughout, $(b_1,b_2)$ is an admissible pair, $b_1=(p,q)$, $b_2=(r,s)$, and $c=(r,q)$.

\begin{definition}[The two slides]
\label{def:slides}
\begin{enumerate}
\item \textbf{The vertical slide of $b_1$:} put $D'=D\setminus\{b_1\}$ and
\[
\pi_V(i,j)=
\begin{cases}
(i+1,j), & \text{if $c\in D'$, $i\ge r$ and $j\ge q$},\\[2pt]
(i,j),   & \text{otherwise};
\end{cases}
\qquad
D_1=\pi_V(D')\cup\{c\}.
\]
That is, $b_1$ is moved to $c$, and if a box already occupies $c$ then every box lying weakly to the lower-right of $c$ is translated one step downward.

\item \textbf{The horizontal slide of $b_2$:} put $D''=D\setminus\{b_2\}$ and
\[
\pi_H(i,j)=
\begin{cases}
(i,j+1), & \text{if $c\in D''$, $i\ge r$ and $j\ge q$},\\[2pt]
(i,j),   & \text{otherwise};
\end{cases}
\qquad
D_2=\pi_H(D'')\cup\{c\}.
\]
That is, $b_2$ is moved to $c$, and if a box already occupies $c$ then every box lying weakly to the lower-right of $c$ is translated one step rightward.
\end{enumerate}
In both cases the induced maps of boxes
\begin{align*}
 \phi_V\colon D\to D_1,\qquad \phi_V(\beta)=\begin{cases}c,&\beta=b_1,\\ \pi_V(\beta),&\beta\in D',\end{cases}\\
 \phi_H\colon D\to D_2,\qquad
 \phi_H(\beta)=\begin{cases}c,&\beta=b_2,\\ \pi_H(\beta),&\beta\in D'',\end{cases}
\end{align*}
are bijections, and the induced maps on tableaux are denoted by $\sigma_V$ and $\sigma_H$:
\[
(\sigma_V T)(\phi_V(\beta))=T(\beta),
\qquad
(\sigma_H T)(\phi_H(\beta))=T(\beta)
\qquad(\beta\in D).
\]
\end{definition}
The repeated application of the sliding operations to a diagram $D$ is naturally organized into a rooted binary tree, which we call the \emph{sliding tree} of $D$ and denote by $\mathcal{T}_D$; for $D=D(\lambda)$ we also write $\mathcal{T}_\lambda$.
The root of $\mathcal{T}_D$ is $D$ itself.
Each internal node of $\mathcal{T}_D$ is a diagram $E$ possessing a pre-ordered box, at which a consecutive admissible pair $(b_1,b_2)$ has been fixed.
Then, the node has exactly two children, namely $E_1$, obtained from $E$ by the vertical slide of $b_1$, and $E_2$, obtained from $E$ by the horizontal slide of $b_2$. 
A leaf of $\mathcal{T}_D$ is a diagram having no pre-ordered box; 
if $D$ is corner-crossed, the tree is finite by Corollary~\ref{cor:termination} below. 

We write $\mathcal{P}_D$, or simply $\mathcal{P}$ when $D$ is clear from the context, for the set of root-to-leaf paths of $\mathcal{T}_D$.
For corner-crossed $D$ and $P\in\mathcal{P}_D$, we denote by $\alpha_P$ the composition of the anti-ribbon diagram at the leaf of $P$ and by
\[
\Phi_P\colon D\longrightarrow D_{\mathrm{ar}}(\alpha_P),
\qquad
\sigma_P\colon T(D,X)\longrightarrow T_{\mathrm{ar}}(\alpha_P,X)
\]
the composites along the edges of $P$ of the box bijections $\phi_V,\phi_H$ and of the induced maps $\sigma_V,\sigma_H$ of Definition~\ref{def:slides}, diagrams being identified with their normalizations.

Throughout the paper we depict one slide, and its iteration, by an equation between labelled diagrams such as the one in the example below, in which the letters stand for weights and not for fixed numbers.

\begin{example}
\label{ex1}
Let $\lambda=(2,1)$ and let
\[
A=
\begin{ytableau}
a & b \\
c
\end{ytableau}
\]
be a tableau of shape $\lambda$.
The pre-ordered boxes of $D(\lambda)$ are $(1,2)$ and $(2,1)$; both are outer, and $(b_1,b_2)=\bigl((1,2),(2,1)\bigr)$ is the unique admissible pair, with
\[
p=1,\quad q=2,\quad r=2,\quad s=1,\qquad (r,q)=(2,2)\notin D(\lambda).
\]
Because the cell $(r,q)$ is unoccupied, neither slide involves a translation: 
the horizontal slide moves $b_2$ to $(2,2)$ and the vertical slide moves $b_1$ to $(2,2)$.
Both resulting diagrams have no pre-ordered box, so the sliding tree has a single internal node and two leaves, and
\begin{align*}
\begin{ytableau}
a & b \\
c
\end{ytableau}
&=
\begin{ytableau}
a & b \\
\none & c
\end{ytableau}
+
\begin{ytableau}
a \\
c & b
\end{ytableau}.
\end{align*}
The two leaves are the anti-ribbon diagrams associated with the compositions $\alpha_1=(2,1)$ and $\alpha_2=(1,2)$, with words $\mathtt{srb}$ and $\mathtt{sbr}$ respectively; the corresponding conditions on the entries are
\[
a\le b<c
\qquad\text{and}\qquad
a<c\le b,
\]
which are mutually exclusive and exhaust all semi-standard fillings of $\lambda$.
In accordance with Theorems~\ref{thm:onetoone} and~\ref{thm:hook}, the number of terms is $f^{(2,1)}=3!/(3\cdot1\cdot1)=2$.
\end{example}

\begin{lemma}
\label{lem:position}
Let $(b_1,b_2)$ be an admissible pair of $D$, $b_1=(p,q)$, $b_2=(r,s)$, so that $p<r$ and $q>s$. Then:
\begin{enumerate}
\item In row $r$, there is no box of $D$ in the open column range $(s,q)$. A box may occur at column $q$ itself.
\item In column $q$, there is no box of $D$ in the open row range $(p,r)$. A box may occur at row $r$ itself.
\end{enumerate}
\end{lemma}

\begin{proof}
We prove the first assertion. Suppose that $(r,t)\in D$ for some $s<t<q$.
Because $p<r$ and $t<q$, we have $(p-r)(t-q)>0$, so that $\{b_1,(r,t)\}$ is a crossing and $(r,t)$ is pre-ordered.
The boxes $b_2$ and $(r,t)$ lie in the row $r$ with $s<t$, so that $b_2\preceq_D(r,t)$, contradicting that $b_2=(r,s)$ is outer pre-ordered.
One may prove the second assertion by a symmetric argument.
\end{proof}

The next lemma shows that the procedure never stalls: as long as a crossing remains, an admissible pair exists; it may always be chosen to consist of two \emph{consecutive} outer pre-ordered boxes, and its two boxes are extremal in their row and in their column, respectively.
\begin{lemma}
\label{lem:admissible}
Let $D$ be a diagram with $\operatorname{Pre}(D)\neq\emptyset$.
\begin{enumerate}
\item There exists an admissible pair $(b_1,b_2)=((p,q),(r,s))$. One may take
\[
q=\max\{\,j\mid (i,j)\in\operatorname{Pre}(D)\,\},
\qquad
p=\max\{\,i\mid (i,q)\in\operatorname{Pre}(D)\,\},
\]
and then
\[
r=\max\{\,i\mid (i,j)\in\operatorname{Pre}(D),\ i>p,\ j<q\,\},
\qquad
s=\max\{\,j\mid (r,j)\in\operatorname{Pre}(D),\ j<q\,\}.
\]
\item Listing the outer pre-ordered boxes as $b^{(1)},\ldots,b^{(m)}$ in increasing order of their rows as in Definition~\ref{def:crossing}, at least one consecutive pair $(b^{(t)},b^{(t+1)})$ is admissible.
\item Let $(b_1,b_2)$ be any admissible pair, $b_1=(p,q)$, $b_2=(r,s)$. Then, $b_1$ is the rightmost box of row $p$ of $D$, and $b_2$ is the lowest box of column $s$ of $D$.
\end{enumerate}
\end{lemma}

\begin{proof}
(1) Because $\operatorname{Pre}(D)\neq\emptyset$ and $D$ is finite, $q$ and $p$ are well defined and $b_1=(p,q)\in\operatorname{Pre}(D)$.
Let $v=(i,j)\in\operatorname{Pre}(D)$ with $b_1\preceq_D v$. Every step of a chain realizing $\preceq_D$ increases one coordinate and leaves the other unchanged, so $i\ge p$ and $j\ge q$; the maximality of $q$ forces $j=q$ and then the maximality of $p$ forces $i=p$. Hence $v=b_1$ and $b_1\in\operatorname{Out}(D)$.

Because $b_1$ is pre-ordered, there is $w=(i',j')\in D$ with $(i'-p)(q-j')>0$. If $i'<p$ and $j'>q$, then $w$ is pre-ordered and lies in a column $j'>q$, contradicting the maximality of $q$. Hence $i'>p$ and $j'<q$, and $w$ is pre-ordered, so the set
\[
S=\{\,(i,j)\in\operatorname{Pre}(D)\mid i>p,\ j<q\,\}
\]
is non-empty; thus $r$ and $s$ are well defined and $b_2=(r,s)\in S$ satisfies $p<r$ and $s<q$.

It remains to check $b_2\in\operatorname{Out}(D)$. Let $v=(i,j)\in\operatorname{Pre}(D)$ with $b_2\preceq_D v$, so that $i\ge r$ and $j\ge s$ as above. If $j<q$, then $i\ge r>p$ shows $v\in S$, and the maximality of $r$ gives $i=r$, after which the maximality of $s$ gives $j=s$, so $v=b_2$. If $j=q$, then $(i,q)\in\operatorname{Pre}(D)$ with $i\ge r>p$ contradicts the maximality of $p$; and $j>q$ contradicts the maximality of $q$.
(2) Write $b^{(k)}=(p_k,q_k)$ with $p_1<\cdots<p_m$.
By (1) some pair $(b^{(k)},b^{(l)})$ with $k<l$ is admissible, so $q_k>q_l$. Consequently $q_t>q_{t+1}$ for at least one $t$ with $k\le t<l$, and for such a $t$ the pair $(b^{(t)},b^{(t+1)})$ satisfies $p_t<p_{t+1}$ and $q_t>q_{t+1}$, hence is admissible.

(3) Suppose $(p,J)\in D$ with $J>q$. By $r>p$ and $s<q<J$, the pair $\{(p,J),b_2\}$ is a crossing, so $(p,J)\in\operatorname{Pre}(D)$, and $b_1\preceq_D(p,J)$ since the two boxes lie in the row $p$; this contradicts $b_1\in\operatorname{Out}(D)$. Suppose $(I,s)\in D$ with $I>r$. Because $I>r>p$ and $s<q$, the pair $\{(I,s),b_1\}$ is a crossing, so $(I,s)\in\operatorname{Pre}(D)$, and $b_2\preceq_D(I,s)$ since the two boxes lie in the column $s$; this contradicts $b_2\in\operatorname{Out}(D)$.
\end{proof}
\begin{lemma}
\label{lem:quadrant}
Let $(b_1,b_2)$ be an admissible pair of $D$, $b_1=(p,q)$, $b_2=(r,s)$, and put $c=(r,q)$. If $c\in D$, then no box of $D$ has row index $>r$ and column index $<q$.
\end{lemma}

\begin{proof}
Such a box $v=(i,j)$ would satisfy $(i-r)(q-j)>0$, so $\{v,c\}$ would be a crossing and $c\in\operatorname{Pre}(D)$; but $b_1$ and $c$ lie in the column $q$ with $p<r$, whence $b_1\preceq_D c$, contradicting $b_1\in\operatorname{Out}(D)$.
\end{proof}

\begin{definition}[Corner-crossed]
\label{def:nwcorner}
A box $x=(i,j)$ of a diagram $D$ is a \emph{northwest corner} of $D$ if no box of $D$ lies in the column $j$ above $x$, and no box of $D$ lies in the row $i$ to the left of $x$. We write $\operatorname{NW}(D)$ for the set of northwest corners; it is never empty, and no two northwest corners share a row or a column.
We say that $D$ is \emph{corner-crossed} if any two of its northwest corners form a crossing, that is, if listing $\operatorname{NW}(D)$ by increasing row index the column indices strictly decrease.
We say that $D$ is \emph{edge-connected} if any two of its boxes are joined by a sequence of boxes of $D$ in which any two consecutive boxes share an edge.
\end{definition}

\begin{lemma}
\label{lem:nw-classes}
Young diagrams, skew Young diagrams, whether connected or not, and edge-connected diagrams are corner-crossed.
\end{lemma}

\begin{proof}
For a Young diagram, $(1,1)$ is the only northwest corner.
For $D(\lambda/\mu)$, a northwest corner of a non-empty row $i$ can only be the leftmost box $(i,\mu_i+1)$ of that row; if the rows $i<i'$ both carry corners, then $\mu_i\ge\mu_{i'}$ because $\mu$ is weakly decreasing, and the two column indices are distinct, so $\mu_i+1>\mu_{i'}+1$.
Let now $D$ be edge-connected and suppose that $x$ and $y=(i',j')$ are northwest corners with $x$ in the region $\{(a,b)\mid a<i',\ b<j'\}$. Choose a chain $x=z_0,z_1,\ldots,z_m=y$ of boxes of $D$ whose consecutive members are edge-adjacent, write $z_l=(a_l,b_l)$, and let $k$ be minimal with $a_k\ge i'$ or $b_k\ge j'$; it exists because $z_m=y$ does, and $k\ge1$. Because $z_{k-1}$ satisfies $a_{k-1}<i'$ and $b_{k-1}<j'$, and because edge-adjacency changes one coordinate by $1$, either $z_k$ lies in the row $i'$ with $b_k<j'$, which is a box of $D$ to the left of $y$, or $z_k$ lies in the column $j'$ with $a_k<i'$, which is a box of $D$ above $y$; both contradict that $y$ is a northwest corner.
Hence no northwest corner lies strictly to the north-west of another, and $D$ is corner-crossed.
\end{proof}

\begin{lemma}
\label{lem:nw-preserved}
Let $D$ be a diagram, let $(b_1,b_2)$ be a 
admissible pair, $b_1=(p,q)$, $b_2=(r,s)$, and let $\varphi$ be the box bijection $\phi_V$ or $\phi_H$
of the corresponding slide. Then
\[
\operatorname{NW}(\varphi(D))
\subseteq\varphi\bigl(\operatorname{NW}(D)\bigr),
\]
and if $D$ is corner-crossed then so is $\varphi(D)$.
\end{lemma}
\begin{proof}
We treat the vertical slide; 
the horizontal one is obtained by transposing, which exchanges the two slides and preserves both the notion of a northwest corner and that of a crossing.
Write $c=(r,q)$ and recall that $\pi_V$ translates the boxes of $Q=\{(a,b)\mid a\ge r,\ b\ge q\}$ one step downwards, and only when $c\in D$.

We first show that no new northwest corner appears. The box $c=\varphi(b_1)$ is not one, because $b_2=(r,s)$ has $s<q$, is not translated and therefore lies in the row $r$ to its left. Let then $x=(i,j)\in D$ with $x\neq b_1$ fail to be a northwest corner of $D$; 
we check that $\varphi(x)$ fails to be one in $\varphi(D)$.

Suppose first that some $u=(a,j)\in D$ satisfies $a<i$. If $u=b_1$, then $j=q$ and $a=p$, and Lemma~\ref{lem:position}(2) gives $i\ge r$; if $c\notin D$ then $i>r$ and $\varphi(x)=x$, while if $c\in D$ then $x\in Q$ and $\varphi(x)=(i+1,q)$. In both cases $c=(r,q)$ lies in the column $q$ above $\varphi(x)$. If $u\neq b_1$, then $u$ and $x$ lie in the same column, so $u\in Q$ implies $x\in Q$, and the translation being rigid, $\varphi(u)$ still lies in the column $j$ above $\varphi(x)$.

Suppose next that some $v=(i,b)\in D$ satisfies $b<j$. By Lemma~\ref{lem:admissible}(3) the box $b_1$ is the rightmost box of its row, so $v\neq b_1$. If $v$ and $x$ are both translated, or both fixed, then $\varphi(v)$ still lies in the row of $\varphi(x)$ to its left. Otherwise $x\in Q$ and $v\notin Q$, so that $c\in D$, $i\ge r$ and $b<q\le j$, whence $i=r$ by Lemma~\ref{lem:quadrant}. If $j=q$ then $x=c$ and $\varphi(b_1)=(r,q)$ lies in the column $q$ above $\varphi(x)=(r+1,q)$. If $j>q$ then $c$ lies in the row $r$ with $q\le j$, so $c\in Q$ is translated together with $x$ and $\varphi(c)$ lies in the row of $\varphi(x)$ to its left.

Now let $D$ be corner-crossed and let $\varphi(x)$, $\varphi(y)$
be two northwest corners of $\varphi(D)$, with $x=(i,j)$ and $y=(i',j')$ northwest corners of $D$, say $i<i'$; then $j>j'$. 
The vertical slide preserves the column of every box, so the columns of $\varphi(x)$ and $\varphi(y)$ are again $j$ and $j'$ with $j>j'$. 
If $x$ and $y$ are both translated, or both fixed, their rows become $i<i'$ or $i+1<i'+1$. If $y$ is translated and $x$ is not, the rows are $i<i'<i'+1$. 
If $x$ is translated and $y$ is not, then $i\ge r$, hence $i'>i\ge r$, and $y\notin Q$ forces $j'<q$; the rows become $i+1$ and $i'$ with $i+1\le i'$, and $i+1=i'$ is impossible, since $\varphi(x)$ and $\varphi(y)$ would then lie in the same row with $j'<j$, so that $\varphi(y)$ would lie to the left of $\varphi(x)$. 
In all cases $\varphi(x)$ and $\varphi(y)$ form a crossing.
\end{proof}

For corner-crossed diagrams, the row--column order is simply the componentwise order.
\begin{lemma}
\label{lem:componentwise}
Let $D$ be a corner-crossed diagram. 
For $u=(a,b)$ and $v=(c,d)$ in $D$, one has $u\preceq_D v$ if and only if $a\le c$ and $b\le d$.
\end{lemma}

\begin{proof}
Every step of a chain realizing $\preceq_D$ increases one coordinate and leaves the other unchanged, which gives the implication from left to right.
Conversely, let $a\le c$ and $b\le d$; we argue by induction on $c-a$. Let $(c,m)$ be the leftmost box of the row $c$ of $D$, so that $m\le d$.
If $m\le b$, then the row $c$ contains $(c,m)$ and $v$ with $m\le b\le d$, so that $(c,b)\in\widetilde{D}$; the column $b$ then contains $u$ and $(c,b)$ with $a\le c$, and $u\preceq_D(c,b)\preceq_D v$.
If $m>b$, then $a<c$, since otherwise $u$ would be a box of the row $c$ to the left of $(c,m)$; as $(c,m)\preceq_D v$ holds in the row $c$, it suffices to show $u\preceq_D(c,m)$.
The box $(c,m)$ is not a northwest corner of $D$: otherwise, starting from $u$ and passing repeatedly to a box of $D$ above in the same column, or to the left in the same row, as long as one exists, we would reach a northwest corner $y$ of $D$ with row index $\le a<c$ and column index $\le b<m$, and $y$ and $(c,m)$ would not form a crossing, contradicting the hypothesis.
Since no box of $D$ lies in the row $c$ to the left of $(c,m)$, there is thus a box $(i,m)\in D$ with $i<c$.
If $i\le a$, then the column $m$ contains $(i,m)$ and $(c,m)$ with $i\le a<c$, so that $(a,m)\in\widetilde{D}$; the row $a$ then contains $u$ and $(a,m)$ with $b<m$, and $u\preceq_D(a,m)\preceq_D(c,m)$.
If $i>a$, then $i-a<c-a$ and $b<m$, so the induction hypothesis gives $u\preceq_D(i,m)$, and $(i,m)\preceq_D(c,m)$ holds in the column $m$.
\end{proof}

\begin{lemma}
\label{lem:gap}
Let $D$ be a corner-crossed diagram and let $(b_1,b_2)$ be an admissible pair, $b_1=(p,q)$, $b_2=(r,s)$.
Then $D$ has no box $(i,j)$ with $p<i<r$ and $j>q$, and no box $(i,j)$ with $i>r$ and $s<j<q$.
\end{lemma}

\begin{proof}
Suppose first that $x=(i,j)\in D$ satisfies $p<i<r$ and $j>q$.
Since $i<r$ and $j>q>s$, the boxes $x$ and $b_2$ form a crossing, and hence $x\in\operatorname{Pre}(D)$.
Moreover, by Lemma \ref{lem:componentwise}, we have
$b_1\preceq_D x$ with $x\neq b_1$.
This contradicts the maximality of $b_1$ in $\operatorname{Pre}(D)$.

Similarly, suppose that $x=(i,j)\in D$ satisfies $i>r$ and $s<j<q$.
Then $i>r>p$ and $j<q$, so $x$ and $b_1$ form a crossing, and hence $x\in\operatorname{Pre}(D)$.
Since $b_2\preceq_D x$ with $x\neq b_2$, this contradicts the maximality of $b_2$ in $\operatorname{Pre}(D)$.
\end{proof}

The following lemma is the key to the construction: a slide destroys the crossing formed by the admissible pair and creates no new one, so that the number of crossings strictly decreases at every step.

\begin{lemma}
\label{lem:no-new-crossing}
Let $D$ be a corner-crossed diagram, let $(b_1,b_2)$ be an admissible pair of $D$, and let $D_1,\phi_V$ and $D_2,\phi_H$ be as in Definition~\ref{def:slides}. Then
\[
\operatorname{Cr}(D_1)\subseteq
\phi_V\bigl(\operatorname{Cr}(D)\setminus\bigl\{\{b_1,b_2\}\bigr\}\bigr),
\qquad
\operatorname{Cr}(D_2)\subseteq
\phi_H\bigl(\operatorname{Cr}(D)\setminus\bigl\{\{b_1,b_2\}\bigr\}\bigr),
\]
where $\phi_V,\phi_H$ act on pairs componentwise. In particular
\[
|\operatorname{Cr}(D_1)|\le|\operatorname{Cr}(D)|-1,
\qquad
|\operatorname{Cr}(D_2)|\le|\operatorname{Cr}(D)|-1.
\]
\end{lemma}

\begin{proof}
We treat the vertical slide. The horizontal case is obtained by transposing. Write $b_1=(p,q)$, $b_2=(r,s)$ and $c=(r,q)$.

We first show that the translation $\pi_V$ creates no crossing.
If $c\notin D'$, then $\pi_V$ is the identity and there is nothing to prove.
Assume $c\in D'$, and let $E_+=\{(i,j)\in D'\mid i\ge r,\ j\ge q\}$ be the translated set. 
Crossings between two boxes of $E_+$, and crossings between two boxes of $D'\setminus E_+$, are unaffected. Let $X=(i,j)\in D'\setminus E_+$ and $Y=(i',j')\in E_+$, and put $Y'=(i'+1,j')$.
If $\{X,Y'\}$ is a crossing with $Y'$ in the third quadrant with respect to $X$, then $i'+1>i$ and $j'<j$; when $i'>i$ the pair $\{X,Y\}$ was already a crossing, and when $i'=i$ we get $j>j'\ge q$ and $i=i'\ge r$, so $X\in E_+$, a contradiction.
If $\{X,Y'\}$ is a crossing with $Y'$ in the first quadrant with respect to $X$, then $i'+1<i$ and $j'>j$, so already $i'<i$ and $j'>j$ and $\{X,Y\}$ was a crossing.
By the definition of a crossing, these are the only two cases.
Hence $\operatorname{Cr}(\pi_V(D'))\subseteq\pi_V(\operatorname{Cr}(D'))$.

We next show that inserting $c$ creates no crossing.
Let $X=(i,j)\in\pi_V(D')$ and suppose $\{X,c\}\in\operatorname{Cr}(D_1)$.
If $i<r$ and $j>q$, then $X$ lies in a row $<r$ and is therefore not translated, 
so $X\in D$. Since $b_1$ is the rightmost box of its row by Lemma~\ref{lem:admissible}(3), we have $i\neq p$, and $p<i<r$ is excluded by Lemma~\ref{lem:gap}; hence $i<p$. Together with $j>q$, this shows that $\{X,b_1\}$ is a crossing of $D$, and $\{X,c\}=\phi_V(\{X,b_1\})$.
If $i>r$ and $j<q$, then $X$ lies in a column $<q$ and is not translated, so $X\in D$. By $i>r>p$ and $j<q$, the pair $\{X,b_1\}$ is a crossing of $D$ and $\{X,c\}=\phi_V(\{X,b_1\})$.

Note that $X\neq b_2$ in either case.
Indeed, $b_2=(r,s)$ with $s<q$ is not translated, and $b_2$ and $c$ lie in the same row $r$, so $\{b_2,c\}$ is not a crossing of $D_1$.
Hence $\{X,b_1\}\in\operatorname{Cr}(D)\setminus\{\{b_1,b_2\}\}$.
Since $b_1\notin D'$, we have $\operatorname{Cr}(D')\subseteq\operatorname{Cr}(D)\setminus\{\{b_1,b_2\}\}$, and the two steps above prove the asserted inclusion.
Finally, $\{b_1,b_2\}\in\operatorname{Cr}(D)$ and $\phi_V$ is injective, so the asserted inequality follows.
\end{proof}

Because every crossing of the new diagram is the image of a crossing of $D$, a box lying in no crossing of $D$ lies in no crossing of any descendant of $D$ in the sliding tree. The number of crossings also provides a termination measure for the procedure.

\begin{corollary}
\label{cor:termination}
Let $D$ be a corner-crossed diagram. 
Starting from $D$ and applying slides with respect to consecutive admissible pairs repeatedly, every branch reaches, after at most $|\operatorname{Cr}(D)|$ steps, a diagram $E$ with $\operatorname{Cr}(E)=\emptyset$. 
Consequently the sliding tree $\mathcal{T}_D$ is a finite rooted binary tree, and its internal nodes are exactly the nodes admitting an admissible pair.
\end{corollary}

\begin{proof}
By Lemma~\ref{lem:nw-preserved} and induction along the edges, every node is corner-crossed.
Hence, by Lemma~\ref{lem:no-new-crossing}, the non-negative integer $|\operatorname{Cr}(\cdot)|$ strictly decreases along every edge, and by Lemma~\ref{lem:admissible}(2) a node $E$ with $\operatorname{Cr}(E)\neq\emptyset$ admits a consecutive admissible pair.
\end{proof}

\begin{remark}
The hypothesis that $D$ is corner-crossed cannot be omitted from Lemma~\ref{lem:no-new-crossing}.
For $D=\{(1,3),(2,1),(3,2)\}$, whose northwest corners $(2,1)$ and $(3,2)$ do not form a crossing, all three boxes are outer pre-ordered, and $((1,3),(2,1))$ is the unique consecutive admissible pair.
Its horizontal slide yields $\{(1,3),(2,3),(3,2)\}$, which again has two crossings.
\end{remark}

\begin{example}
\label{ex:afull}
We illustrate the whole procedure on the diagram
\[
D=\{(1,1),(1,2),(1,3),(2,1),(2,3)\},
\]
whose boxes we label $a,\ldots,e$ in the order $(1,1),(1,2),(1,3),(2,1),(2,3)$. 
The diagram $D$ is edge-connected, hence corner-crossed by Lemma~\ref{lem:nw-classes}, and it is not a skew Young diagram, since the cell $(2,2)$ is an implicit box.

Its pre-ordered boxes are $(1,2)$, $(1,3)$ and $(2,1)$, of which $(1,3)$ and $(2,1)$ are outer; they form the unique admissible pair, and the cell $c=(2,3)$ is occupied by the box labelled $e$.
Hence the vertical slide of $b_1=(1,3)$ translates $e$ one step downward, and the horizontal slide of $b_2=(2,1)$ translates it one step rightward.
The second resulting diagram has no pre-ordered box.
In the first one, $(1,2)$ and $(2,1)$ form the unique admissible pair and the cell $c=(2,2)$ is unoccupied, so one more round of slides gives
\begin{align*}
\begin{ytableau}
a & b & c \\
d & \none & e
\end{ytableau}
&=
\begin{ytableau}
a & b \\
d & \none & c \\
\none & \none & e
\end{ytableau}
+
\begin{ytableau}
a & b & c \\
\none & \none & d & e
\end{ytableau}\\
&=
\begin{ytableau}
a \\
d & b & c \\
\none & \none & e
\end{ytableau}
+
\begin{ytableau}
a & b \\
\none & d & c \\
\none & \none & e
\end{ytableau}
+
\begin{ytableau}
a & b & c \\
\none & \none & d & e
\end{ytableau}.
\end{align*}
The three leaves are the anti-ribbon diagrams with words $\mathtt{sbrrb}$, $\mathtt{srbrb}$ and $\mathtt{srrbr}$, in agreement with $|\mathrm{SYT}(D)|=3$.
\end{example}

\begin{definition}[Normalization]
\label{def:normalization}
Let $E$ be a diagram. For $\beta\in E$, let $\rho(\beta)$ denote the rank of the row index of $\beta$ among the row indices occurring in $E$, and let $\gamma(\beta)$ denote the rank of the column index of $\beta$ among the column indices occurring in $E$. The \emph{normalization} of $E$ is
\[
\widehat{E}=\{(\rho(\beta),\gamma(\beta))\mid\beta\in E\}.
\]
Normalization preserves the relative order of rows and of columns, hence preserves $\operatorname{Cr}$, $\operatorname{Pre}$ and the order $\preceq$ of Definition~\ref{def:implicit}.
\end{definition}

We now analyze the diagrams occurring at the leaves, namely those with no crossing at all. Such a diagram is totally ordered by the componentwise order, and each pair of consecutive boxes is related in one of three ways; the normalization is an anti-ribbon diagram exactly when the third, diagonal, case never occurs.

\begin{lemma}
\label{lem:leaf-structure}
Let $E$ be a diagram with $\operatorname{Cr}(E)=\emptyset$ and $|E|=N$. 
Then:
\begin{enumerate}
\item the componentwise order, defined by $(i,j)\le(i',j')$ if and only if $i\le i'$ and $j\le j'$, is a total order on $E$; write $E=\{u_1<u_2<\cdots<u_N\}$;
\item for each $k$ exactly one of the following holds: \textup{(r)} $u_k$ and $u_{k+1}$ lie in the same row; \textup{(b)} they lie in the same column; \textup{(d)} both coordinates increase strictly;
\item in case \textup{(r)} no box of $E$ lies in a column strictly between those of $u_k$ and $u_{k+1}$, in case \textup{(b)} no box of $E$ lies in a row strictly between those of $u_k$ and $u_{k+1}$, and in case \textup{(d)} no box of $E$ lies in a row, nor in a column, strictly between those of $u_k$ and $u_{k+1}$;
\item consequently $\widehat{E}$ is an anti-ribbon diagram in the southeast convention if and only if case \textup{(d)} occurs for no $k$, and in that case the associated word is $s\varepsilon_1\cdots\varepsilon_{N-1}$ with $\varepsilon_k=r$ in case \textup{(r)} and $\varepsilon_k=b$ in case \textup{(b)}.
\end{enumerate}
\end{lemma}

\begin{proof}
(1) Let $u=(i,j)$ and $v=(i',j')$ be distinct boxes of $E$. Note that $\{u,v\}$ is not a crossing. Hence, neither ($i<i'$ and $j>j'$) nor ($i>i'$ and $j<j'$) holds, so that $i\le i'$ and $j\le j'$, or $i\ge i'$ and $j\ge j'$. Thus any two boxes are comparable, and $\le$ is antisymmetric because a box is determined by its coordinates.

(2) Clear from (1), the three cases being $i=i'$, $j=j'$ and $i<i'$, $j<j'$; they are mutually exclusive because $u_k\neq u_{k+1}$.

(3) Assume case (r), say $u_k=(i,j)$ and $u_{k+1}=(i,j'')$ with $j<j''$, and suppose $(a,b)\in E$ with $j<b<j''$. 
If $a<i$ then $\{(a,b),u_k\}$ is a crossing; if $a>i$ then $\{(a,b),u_{k+1}\}$ is a crossing. Both are excluded, so $a=i$; 
but then $u_k<(a,b)<u_{k+1}$, contradicting that $u_k$ and $u_{k+1}$ are consecutive. Case (b) is symmetric.
In case (d), say $u_k=(i,j)$ and $u_{k+1}=(i',j')$ with $i<i'$ and $j<j'$, suppose $(a,b)\in E$ has $i<a<i'$. By (1) the box $(a,b)$ is comparable with $u_k$ and with $u_{k+1}$, and these being consecutive we have $(a,b)\le u_k$ or $(a,b)\ge u_{k+1}$, forcing $a\le i$ or $a\ge i'$; both are excluded. The assertion about columns is proved in the same way.

(4) By (3) and Definition~\ref{def:normalization}, in case (r) the boxes $\widehat{u_k}$ and $\widehat{u_{k+1}}$ are horizontally adjacent, and in case (b) vertically adjacent. If case (d) never occurs, $\widehat{u_1},\ldots,\widehat{u_N}$ therefore form a path of unit steps to the right and downward starting at $\widehat{u_1}=(1,1)$, which is exactly an anti-ribbon diagram, with the stated word. Conversely, an anti-ribbon diagram has no step of type (d).
\end{proof}

\begin{corollary}
\label{cor:leaf-antiribbon}
Let $D$ be corner-crossed in the sense of Definition~\ref{def:nwcorner}, for instance a Young diagram, a skew Young diagram, or an edge-connected diagram. 
Then case \textup{(d)} of Lemma~\ref{lem:leaf-structure}(2) occurs at no leaf of $\mathcal{T}_D$, and the normalization of every leaf of $\mathcal{T}_D$ is a single anti-ribbon diagram.
\end{corollary}

\begin{proof}
By Lemma~\ref{lem:nw-preserved} and induction along the edges, every node of $\mathcal{T}_D$ is corner-crossed.
Let $E=\{u_1<\cdots<u_N\}$ be a leaf. The box $u_1$ is a northwest corner of $E$, since every other box has both indices at least as large and no two boxes share both. Suppose case \textup{(d)} occurred at some $k$, that is, $u_k=(i,j)$ and $u_{k+1}=(i',j')$ with $i<i'$ and $j<j'$. Every $u_l$ with $l\le k$ has both indices at most those of $u_k$, hence strictly smaller than those of $u_{k+1}$, and every $u_l$ with $l\ge k+2$ has both indices at least those of $u_{k+1}$; consequently no box of $E$ lies in the column $j'$ above $u_{k+1}$ or in the row $i'$ to its left, and $u_{k+1}$ is a northwest corner. But $u_1$ and $u_{k+1}$ then have strictly increasing row indices and strictly increasing column indices, so they do not form a crossing, contradicting that $E$ is corner-crossed.
Hence no diagonal step occurs and Lemma~\ref{lem:leaf-structure}(4) applies.
\end{proof}

\subsection{Decomposition of tableaux}
One step of the procedure splits the standard tableaux on a node according to the comparison of the two entries occupying the admissible pair.
The following proposition says, more precisely, that a slide refines the row--column order by exactly one relation. Write $\preceq_D^{V}$ for the transitive closure of $\preceq_D\cup\{(b_2,b_1)\}$ and $\preceq_D^{H}$ for that of $\preceq_D\cup\{(b_1,b_2)\}$.

\begin{proposition}
\label{prop:step}
Let $D$ be a corner-crossed diagram in the sense of Definition~\ref{def:nwcorner} and let $(b_1,b_2)$ be a \emph{consecutive} admissible pair, $b_1=(p,q)$, $b_2=(r,s)$, $p<r$, $q>s$.
Let $D_1$ and $D_2$ be the diagrams obtained from $D$ by the vertical slide of $b_1$ and the horizontal slide of $b_2$, respectively.
Then the box bijections are isomorphisms of ordered sets
\[
\phi_V\colon(D,\preceq_D^{V})\ \xrightarrow{\ \sim\ }\ (D_1,\preceq_{D_1}),
\qquad
\phi_H\colon(D,\preceq_D^{H})\ \xrightarrow{\ \sim\ }\ (D_2,\preceq_{D_2}).
\]
Consequently, for $T\in\mathrm{SYT}(D)$ write $x=T(p,q)$, $y=T(r,s)$, and set
\[
\mathrm{SYT}(D)_{>}=\{T : x>y\},\qquad \mathrm{SYT}(D)_{<}=\{T : x<y\}.
\]
Then, $\sigma_V$ restricts to a bijection $\mathrm{SYT}(D)_{>}\xrightarrow{\sim}\mathrm{SYT}(D_1)$, and $\sigma_H$ restricts to a bijection $\mathrm{SYT}(D)_{<}\xrightarrow{\sim}\mathrm{SYT}(D_2)$.
\end{proposition}

\begin{proof}
The last assertion follows from the first. The boxes $b_1$ and $b_2$ form a crossing, hence are $\preceq_D$-incomparable, so a tableau $T\in\mathrm{SYT}(D)$ lies in $\mathrm{SYT}(D)_{>}$ if and only if it is a linear extension of $\preceq_D$ in which $b_2$ precedes $b_1$, that is, a linear extension of $\preceq_D^{V}$; transporting entries along $\phi_V$, which is what $\sigma_V$ does, identifies $\mathrm{SYT}(D)_{>}$ with $\mathrm{SYT}(D_1)$, and the horizontal case is the same.

We prove the first assertion for the vertical slide.
The horizontal one follows by transposing, which exchanges the two slides and preserves corner-crossed diagrams, the order of Lemma~\ref{lem:componentwise} and admissible pairs.
The transposed pair is again consecutive: an outer pre-ordered box $v=(g,h)$ with $s<h<q$ would satisfy $g<p$ or $g>r$ by consecutiveness, and then $v\preceq_D b_1$ or $b_2\preceq_D v$ by Lemma~\ref{lem:componentwise}, contradicting that $v$, $b_1$ and $b_2$ are outer.

Write $c=(r,q)$ and $Q=\{(a,b)\mid a\ge r,\ b\ge q\}$.
The map $\phi_V$ preserves the column of every box; it sends $b_1$ to $c$, translates the boxes of $D\cap Q$ one row downwards when $c\in D$, and fixes all other boxes.
By Lemma~\ref{lem:nw-preserved}, $D_1$ is corner-crossed, so by Lemma~\ref{lem:componentwise} both $\preceq_D$ and $\preceq_{D_1}$ are the componentwise order $\le$.
Since $b_1\le y\le b_2$ is impossible because $q>s$, the relation $\preceq_D^{V}$ is given by
\[
x\preceq_D^{V}y
\iff
x\le y\ \text{ or }\ \bigl(x\le b_2\ \text{ and }\ b_1\le y\bigr).
\]
We use two observations.
First, if $y=(e,f)\in D$ satisfies $b_1\le y$ and $y\neq b_1$, then $y\in Q$. Indeed, $e>p$ because $b_1$ is the rightmost box of its row (Lemma~\ref{lem:admissible}(3)), and $p<e<r$ is excluded by Lemma~\ref{lem:position}(2) if $f=q$ and by Lemma~\ref{lem:gap} if $f>q$; hence $e\ge r$.
Second, if $x\le b_2$, then $x$ lies in a row $\le r$ and in a column $\le s<q$, so that $x\notin Q$ and $\phi_V(x)=x$.
We now show that $\phi_V(x)\le\phi_V(y)$ if and only if $x\preceq_D^{V}y$, for distinct $x,y\in D$.

Suppose first that $x\neq b_1$ and $y\neq b_1$. If $x\le b_2$ and $b_1\le y$, then the two observations give $x\le y$; hence $x\preceq_D^{V}y$ if and only if $x\le y$.
If $c\notin D$, then $\phi_V$ fixes $x$ and $y$.
If $c\in D$ and $x,y$ both lie in $Q$ or both lie outside $Q$, their rows are shifted by the same amount.
If $x\in Q$ and $y\notin Q$, then neither $x\le y$ nor $\phi_V(x)\le\phi_V(y)$ holds, since either inequality would force $y\in Q$.
If $x=(i,j)\notin Q$ and $y=(e,f)\in Q$, then $x\le y$ implies $\phi_V(x)=x\le(e+1,f)=\phi_V(y)$; conversely, $\phi_V(x)\le\phi_V(y)$ means $i\le e+1$ and $j\le f$, and $i=e+1$ is impossible, since then $i>r$ and $j<q$, contrary to Lemma~\ref{lem:quadrant}.
In all cases $\phi_V(x)\le\phi_V(y)$ if and only if $x\le y$.

Suppose next that $x=b_1$. Then $b_1\preceq_D^{V}y$ if and only if $b_1\le y$, because $b_1\le b_2$ fails.
If $b_1\le y$, then $y\in Q$ by the first observation, so $\phi_V(y)$ lies in a row $\ge r$ and in a column $\ge q$, that is, $c\le\phi_V(y)$.
Conversely, if $c\le\phi_V(y)$, then either $y\in Q$ or $y$ is fixed and lies in a row $\ge r$; in both cases $y$ lies in a row $\ge r>p$ and in a column $\ge q$, so that $b_1\le y$.

Suppose finally that $y=b_1$. Then $x\preceq_D^{V}b_1$ if and only if $x\le b_1$ or $x\le b_2$, and in either case $x$ is fixed and $\phi_V(x)=x\le c$.
Conversely, let $\phi_V(x)\le c=(r,q)$. Then $x$ is not translated, so $x=(i,j)$ with $i\le r$ and $j\le q$.
If $i\le p$, then $x\le b_1$.
If $i=r$, then $j\neq q$, since $c\in D$ would be translated, and $s<j<q$ is excluded by Lemma~\ref{lem:position}(1); hence $j\le s$ and $x\le b_2$.
Let $p<i<r$. Then $j\neq q$ by Lemma~\ref{lem:position}(2). Assume $s<j<q$. Then $x$ forms a crossing with $b_1$, hence is pre-ordered, and there is $v=(g,h)\in\operatorname{Out}(D)$ with $x\preceq_D v$, so that $g\ge i>p$ and $h\ge j>s$.
Since $(b_1,b_2)$ is consecutive, no outer pre-ordered box lies in a row strictly between $p$ and $r$, and $v\neq b_2$ because $h>s$; since outer pre-ordered boxes lie in distinct rows, $g>r$.
Then $b_2\le v$ and $b_2\preceq_D v$ with $v\neq b_2$, contradicting $b_2\in\operatorname{Out}(D)$. Therefore $j\le s$ and $x\le b_2$.

Hence $\phi_V$ is an isomorphism of ordered sets from $(D,\preceq_D^{V})$ onto $(D_1,\preceq_{D_1})$.
\end{proof}

\begin{theorem}
\label{thm:onetoone}
Let $D$ be a corner-crossed diagram with $n$ boxes, for instance the Young diagram of a partition of $n$, a skew Young diagram or an edge-connected diagram.
There is a bijection between $\mathrm{SYT}(D)$ and the set of root-to-leaf paths in $\mathcal{T}_D$ such that, for each $T\in\mathrm{SYT}(D)$, transporting the entries of $T$ along the corresponding sequence of slides produces the entries $1,2,\ldots,n$ in order along the anti-ribbon path at the leaf.
\end{theorem}

\begin{proof}
Let $k=|\operatorname{Cr}(D)|$ denote the number of crossings of $D$; we proceed by induction on $k$.
If $k=0$, then $D$ has no pre-ordered box, so that $\mathcal{T}_D$ consists of the root $D$ alone, and the normalization of $D$ is an anti-ribbon diagram by Corollary~\ref{cor:leaf-antiribbon}.
By Lemma~\ref{lem:leaf-structure} (1) and Lemma~\ref{lem:componentwise}, the order $\preceq_D$ is total. 
Hence $\mathrm{SYT}(D)$ consists of the single tableau numbering the boxes of $D$ along the anti-ribbon path, and the statement holds.

Suppose $k>0$. 
Let $(b_1,b_2)$ be the consecutive admissible pair fixed at the root of $\mathcal{T}_D$, which exists by Lemma~\ref{lem:admissible}(2), and let $D_1$ and $D_2$ be the diagrams obtained from $D$ by the vertical slide of $b_1$ and the horizontal slide of $b_2$, respectively. 
These are the two children of the root, and the subtrees of $\mathcal{T}_D$ rooted at them are sliding trees of $D_1$ and $D_2$.
By Lemma~\ref{lem:no-new-crossing} we have $|\operatorname{Cr}(D_1)|\le k-1$ and $|\operatorname{Cr}(D_2)|\le k-1$, and by Lemma~\ref{lem:nw-preserved} both diagrams are again corner-crossed, so that the induction hypothesis applies to them. 
By Proposition~\ref{prop:step},
\[
\mathrm{SYT}(D)=\mathrm{SYT}(D)_{>}\,\sqcup\,\mathrm{SYT}(D)_{<}
\]
and this splitting corresponds bijectively, via $\sigma_V$ and $\sigma_H$, to $\mathrm{SYT}(D_1)\sqcup\mathrm{SYT}(D_2)$.
Composing with the bijections furnished by the induction hypothesis and prepending the first edge of $\mathcal{T}_D$, we obtain the desired bijection between $\mathrm{SYT}(D)$ and root-to-leaf paths in $\mathcal{T}_D$, compatible with the anti-ribbon-reading property.
\end{proof}

Having established termination and the shape of the leaves for corner-crossed diagrams, we now turn to semi-standard tableaux.
We do not attach a set of semi-standard tableaux to an intermediate node of the sliding tree.
A standard tableau only records the order in which the boxes are filled, and this order is carried along the slides; a semi-standard tableau also distinguishes weak from strict inequalities, and after sliding, two boxes in the same relative position may carry either one, according to their origin.

\section{Explicit anti-ribbon expansion via descents}
\label{subsec:anti-ribbon-descents}
In this section, we describe an explicit combinatorial form of the anti-ribbon expansion.
The essential observation is that the anti-ribbon shape attached to a standard tableau is determined by its descent set.
This structural fact explains why multiplicities appear when the expansion is grouped by anti-ribbon words.
Anti-ribbon diagrams are understood in the southeast convention of Section~\ref{subsec:antiribbon}.
\begin{definition}
Let $D$ be a diagram with $N$ boxes and let $T\in \operatorname{SYT}(D)$.
For $1\leq k\leq N$, let $\operatorname{row}_T(k)$ be the row index of the box $T^{-1}(k)$.
The descent set of $T$ is defined by
\[
\operatorname{Des}(T)
=
\{\,k\in\{1,\ldots,N-1\}
\mid
\operatorname{row}_T(k+1)>\operatorname{row}_T(k)\,\}.
\]
Equivalently, $k\in \operatorname{Des}(T)$ if and only if the box containing $k+1$ lies in a row strictly below that of the box containing $k$.

The anti-ribbon word, or $rb$-word, associated with $T$ is defined by
\[
\operatorname{rb}(T)
=
s\varepsilon_1\varepsilon_2\cdots \varepsilon_{N-1},
\]
where
\[
\varepsilon_k
=
\begin{cases}
b,& k\in \operatorname{Des}(T),\\
r,& k\notin \operatorname{Des}(T).
\end{cases}
\]
Thus the positions of the letters $b$ in $\operatorname{rb}(T)$ correspond precisely to the descents of $T$.
In the diagram associated with the word, the letter $r$ represents a step to the right and the letter $b$ represents a step downward.
\end{definition}

For $w\in\mathcal{W}_N$ we use the anti-ribbon diagram $R(w)$ of Definition~\ref{def:Rw}, with its boxes $p_1,\ldots,p_N$ in path order; by Lemma~\ref{lem:dictionary}, $R(w)=D_{\mathrm{ar}}(\alpha(w))$ up to translation.

Let $A:D\to \mathcal{A}$ be a tableau over an alphabet $\mathcal{A}$.
For $T\in \operatorname{SYT}(D)$, write
\[
\operatorname{rb}(T)=w.
\]
We define a tableau $R_T(A)$ of shape $R(w)$ by
\[
R_T(A)(p_k)=A(T^{-1}(k)),
\qquad 1\leq k\leq N.
\]
That is, the entry of $A$ in the box numbered $k$ by $T$ is assigned to the $k$-th box of the anti-ribbon path.

\begin{definition}
The anti-ribbon expansion of $A$ is the formal sum
\[
\operatorname{ARE}_D(A)
=
\sum_{T\in \operatorname{SYT}(D)} R_T(A).
\]
This is an equality in the free $\mathbb{Q}$-vector space generated by anti-ribbon tableaux.
It is not an equality of ordinary tableaux.
\end{definition}

The semi-standard fillings are distributed over the terms of this expansion by the following classical standardization argument; see \cite{Gessel1984} and \cite[Section~7.19]{Stanley1999}.

\begin{theorem}
\label{thm:standardization}
Let $D$ be a corner-crossed diagram with $N$ boxes.
For $T\in\mathrm{SYT}(D)$ put $\alpha_T=\alpha(\operatorname{rb}(T))$ and
\[
\mathcal{S}_T
=
\left\{
M\in T(D)\ \middle|\
\begin{aligned}
&M(T^{-1}(k))\le M(T^{-1}(k+1))
&& (1\le k<N),\\
&M(T^{-1}(k))<M(T^{-1}(k+1))
&& (k\in\operatorname{Des}(T))
\end{aligned}
\right\}.
\]
Then the following hold.
\begin{enumerate}
\item $\mathrm{SSYT}(D)=\coprod_{T\in\mathrm{SYT}(D)}\mathcal{S}_T$.
\item For every $T\in\mathrm{SYT}(D)$, the map $M\mapsto R_T(M)$ is a bijection from $\mathcal{S}_T$ onto $\mathrm{SSYT}_{\mathrm{ar}}(\alpha_T)$ which preserves the multiset of entries.
\end{enumerate}
\end{theorem}

\begin{proof}
(1) Fix $T$ and write $u_k=T^{-1}(k)$. We first show $\mathcal{S}_T\subseteq\mathrm{SSYT}(D)$. Let $M\in\mathcal{S}_T$. Let $u=(i,j)$ and $v=(i',j')$ be boxes of $D$ with $u\preceq_D v$ and $u\neq v$, and put $a=T(u)<b=T(v)$. Then $M(u)=M(u_a)\le\cdots\le M(u_b)=M(v)$. Let moreover $i<i'$. If no $k$ with $a\le k<b$ belonged to $\operatorname{Des}(T)$, the row indices of $u_a,u_{a+1},\ldots,u_b$ would be weakly decreasing, which contradicts that $u_a$ lies in the row $i$ and $u_b$ in the row $i'$. Hence some such $k$ is a descent, and $M(u)=M(u_a)\le\cdots\le M(u_k)<M(u_{k+1})\le\cdots\le M(u_b)=M(v)$.

Next, let $M\in\mathrm{SSYT}(D)$.
By Lemma~\ref{lem:componentwise}, $\preceq_D$ is the componentwise order, so two comparable boxes carrying equal entries of $M$ lie in the same row.
Boxes carrying equal entries of $M$ lie in distinct columns, and we let $\operatorname{std}(M)$ be the numbering of $D$ by $1,\ldots,N$ which lists the boxes in increasing order of their entries, boxes with equal entries being listed in increasing order of their columns. 
If $u\preceq_D v$ and $u\neq v$, then either $M(u)<M(v)$, or $M(u)=M(v)$ and $v$ lies to the right of $u$ in the same row; hence $\operatorname{std}(M)\in\mathrm{SYT}(D)$. 
By construction the entries of $M$ are weakly increasing along $\operatorname{std}(M)$. 
Let $k\in\operatorname{Des}(\operatorname{std}(M))$, and let $x=(i,j)$ and $y=(i',j')$ be the boxes numbered $k$ and $k+1$, so that $i<i'$. If $M(x)=M(y)$, then $j<j'$ by the choice of the order, so that $x\preceq_D y$ by Lemma~\ref{lem:componentwise}, and $i<i'$ gives $M(x)<M(y)$, a contradiction. 
Thus $M\in\mathcal{S}_{\operatorname{std}(M)}$.

Finally, let $M\in\mathcal{S}_T$ and suppose $M(u_k)=M(u_{k+1})$, with $u_k=(i,j)$ and $u_{k+1}=(i',j')$. 
Then $k\notin\operatorname{Des}(T)$, so $i'\le i$. If $j'\le j$, then $u_{k+1}\preceq_D u_k$ by Lemma~\ref{lem:componentwise}, and this gives $T(u_{k+1})<T(u_k)$, a contradiction. 
Hence $j<j'$, so that $T$ lists boxes with equal entries in increasing order of their columns, and $T=\operatorname{std}(M)$. 
Therefore the sets $\mathcal{S}_T$ are pairwise disjoint and cover $\mathrm{SSYT}(D)$.

(2) The map $M\mapsto R_T(M)$ is a bijection from $T(D)$ onto the set of fillings of $R(\operatorname{rb}(T))=D_{\mathrm{ar}}(\alpha_T)$, and the entry of $R_T(M)$ at $p_k$ is $M(T^{-1}(k))$. By Lemma~\ref{lem:dictionary}, $S(\alpha_T)=\operatorname{Des}(\operatorname{rb}(T))$, and this set equals $\operatorname{Des}(T)$ by the definition of $\operatorname{rb}(T)$. 
Hence, by Lemma~\ref{lem:ar-ssyt}, $R_T(M)$ is semi-standard if and only if $M\in\mathcal{S}_T$. 
Since $R_T$ only relocates entries, it preserves their multiset.
\end{proof}

\begin{proposition}
\label{prop:grouped-rb}
Let
\[
\mathcal{W}_N
=
\{\,s\varepsilon_1\cdots\varepsilon_{N-1}
\mid
\varepsilon_k\in\{r,b\}\,\}.
\]
For $w\in\mathcal{W}_N$, set
\[
\operatorname{SYT}_w(D)
=
\{\,T\in \operatorname{SYT}(D)\mid \operatorname{rb}(T)=w\,\}.
\]
Then
\[
\operatorname{ARE}_D(A)
=\sum_{w\in\mathcal{W}_N}
\sum_{T\in \operatorname{SYT}_w(D)}
R_T(A).
\]
In particular, under the all-one specialization $A\equiv 1$, the coefficient of the anti-ribbon type $w$ is
\[
c_D(w)
=
|\operatorname{SYT}_w(D)|
=\#\{\,T\in \operatorname{SYT}(D)\mid \operatorname{Des}(T)=\operatorname{Des}(w)\,\}.
\]
\end{proposition}

\begin{proof}
The first identity follows from the disjoint decomposition
\[
\operatorname{SYT}(D)
=
\coprod_{w\in\mathcal{W}_N}\operatorname{SYT}_w(D).
\]
The second assertion follows from the definition of $\operatorname{rb}(T)$, since the positions of the letters $b$ are precisely the descents of $T$.
\end{proof}
For a partition $\lambda$ with conjugate partition $\lambda'$ and a cell $(i,j)\in D(\lambda)$, define its hook length by
\[
h((i,j))
=\lambda_i-j+\lambda'_j-i+1.
\]

\begin{theorem}[Frame--Robinson--Thrall~\cite{FrameRobinsonThrall1954}]
\label{thm:hook}
Let $\lambda$ be a partition of $N$. Then the number $f^\lambda$ of standard Young tableaux of shape $\lambda$ is
\[
f^\lambda
=\frac{N!}{\prod_{(i,j)\in D(\lambda)} h((i,j))}.
\]
\end{theorem}
\begin{corollary}
\label{cor:ungrouped-hook}
Let $\lambda$ be a partition of $N$. 
Then, we have
\[
\sum_{w\in\mathcal{W}_N} c_{D(\lambda)}(w)
=f^\lambda=
\frac{N!}{\prod_{(i,j)\in D(\lambda)} h((i,j))}.
\]
\end{corollary}

\begin{example}
Let $\lambda=(3,3)$.
The standard tableaux of shape $\lambda$ are
\[
\begin{array}{ccccc}
\begin{ytableau}
1&3&5\\
2&4&6
\end{ytableau}
&
\begin{ytableau}
1&2&5\\
3&4&6
\end{ytableau}
&
\begin{ytableau}
1&3&4\\
2&5&6
\end{ytableau}
&
\begin{ytableau}
1&2&4\\
3&5&6
\end{ytableau}
&
\begin{ytableau}
1&2&3\\
4&5&6
\end{ytableau}
\end{array}
\]
and their anti-ribbon words are, respectively,
\[
\mathtt{sbrbrb},
\qquad
\mathtt{srbrrb},
\qquad
\mathtt{sbrrbr},
\qquad
\mathtt{srbrbr},
\qquad
\mathtt{srrbrr}.
\]
Thus, for
\[
A=
\begin{ytableau}
a&b&c\\
d&e&f
\end{ytableau},
\]
the formal expansion is
\[
\begin{ytableau}
a&b&c\\
d&e&f
\end{ytableau}
=
\begin{ytableau}
a\\
d&b\\
\none&e&c\\
\none&\none&f
\end{ytableau}
+
\begin{ytableau}
a&b\\
\none&d&e&c\\
\none&\none&\none&f
\end{ytableau}
+
\begin{ytableau}
a\\
d&b&c\\
\none&\none&e&f
\end{ytableau}
+
\begin{ytableau}
a&b\\
\none&d&c\\
\none&\none&e&f
\end{ytableau}
+
\begin{ytableau}
a&b&c\\
\none&\none&d&e&f
\end{ytableau}.
\]
In this example every anti-ribbon word occurs with multiplicity $1$.
\end{example}

\begin{example}
\label{ex:333}
Let $\lambda=(3,3,3)$.
The hook lengths are
\[
\begin{ytableau}
5&4&3\\
4&3&2\\
3&2&1
\end{ytableau},
\]
so by the hook-length formula
\[
f^{(3,3,3)}
=
\frac{9!}{5\cdot 4\cdot 3\cdot 4\cdot 3\cdot 2\cdot 3\cdot 2\cdot 1}
=
42.
\]
A direct enumeration of standard tableaux of shape $(3,3,3)$ yields a grouped expansion with $40$ distinct anti-ribbon words.
All coefficients equal $1$ except for the two words
\[
\mathtt{srbrbrbrb},
\qquad
\mathtt{sbrbrbrbr},
\]
which have coefficient $2$.
Explicitly, the words with coefficient $1$ are
\[
\begin{array}{llll}
\mathtt{srrbrrbrr},&\mathtt{srrbrbrrb},&\mathtt{srrbrbrbr},&\mathtt{srrbbrrbr},\\
\mathtt{srrbbrbrb},&\mathtt{srbrrbrrb},&\mathtt{srbrrbrbr},&\mathtt{srbrrbbrr},\\
\mathtt{srbrbrrbr},&\mathtt{srbrbrrbb},&\mathtt{srbrbrbrr},&\mathtt{srbrbrbbr},\\
\mathtt{srbrbbrrb},&\mathtt{srbrbbrbr},&\mathtt{srbbrrbrb},&\mathtt{srbbrrbbr},\\
\mathtt{srbbrbrbr},&\mathtt{srbbrbrbb},&\mathtt{srbbrbbrb},&\mathtt{sbrrbrrbr},\\
\mathtt{sbrrbrbrr},&\mathtt{sbrrbrbrb},&\mathtt{sbrrbbrrb},&\mathtt{sbrrbbrbr},\\
\mathtt{sbrbrrbrb},&\mathtt{sbrbrrbbr},&\mathtt{sbrbrbrrb},&\mathtt{sbrbrbrbb},\\
\mathtt{sbrbrbbrr},&\mathtt{sbrbrbbrb},&\mathtt{sbrbbrrbb},&\mathtt{sbrbbrbrb},\\
\mathtt{sbrbbrbbr},&\mathtt{sbbrrbrbr},&\mathtt{sbbrrbbrb},&\mathtt{sbbrbrbrb},\\
\mathtt{sbbrbrbbr},&\mathtt{sbbrbbrbb}.
\end{array}
\]
Hence the sum of the coefficients is
\[
38\cdot 1+2\cdot 2=42=f^{(3,3,3)},
\]
in agreement with Corollary~\ref{cor:ungrouped-hook} and with the interpretation of the coefficient of $w$ as
\[
c_{D(\lambda)}(w)
=\#\{\,T\in \operatorname{SYT}(\lambda)\mid \operatorname{rb}(T)=w\,\}.
\]

We exhibit the two multiplicities explicitly.
The word
\[
\mathtt{srbrbrbrb}
\]
corresponds to the descent set
\[
\{2,4,6,8\},
\]
and the following two standard tableaux of shape $(3,3,3)$ have this descent set:
\[
\begin{array}{cc}
\begin{ytableau}
1&2&4\\
3&6&8\\
5&7&9
\end{ytableau}
&
\begin{ytableau}
1&2&6\\
3&4&8\\
5&7&9
\end{ytableau}
\end{array}.
\]
Similarly, the word
\[
\mathtt{sbrbrbrbr}
\]
corresponds to the descent set
\[
\{1,3,5,7\},
\]
and the following two standard tableaux have this descent set:
\[
\begin{array}{cc}
\begin{ytableau}
1&3&5\\
2&4&7\\
6&8&9
\end{ytableau}
&
\begin{ytableau}
1&3&5\\
2&6&7\\
4&8&9
\end{ytableau}
\end{array}.
\]
This shows that the map
\[
T\longmapsto \operatorname{rb}(T)
\]
from standard tableaux to anti-ribbon words is not injective in general.
Equivalently, distinct standard tableaux may share the same descent set.
\end{example}

\begin{remark}
The example above explains why the appropriate indexing set for the ungrouped expansion is $\operatorname{SYT}(D)$ rather than the set of anti-ribbon words.
The anti-ribbon words index the grouped expansion, and their coefficients are the descent multiplicities.
\end{remark}

\begin{remark}
The expansion is naturally understood on two levels.
The ungrouped level is indexed by $T\in\operatorname{SYT}(D)$, each contributing a single anti-ribbon tableau $R_T(A)$.
The grouped level is indexed by anti-ribbon words, the coefficient of a word $w$ being the number of linear extensions whose anti-ribbon word is $w$.
Multiplicities arise precisely when distinct linear extensions share the same descent set.
\end{remark}

In the following, we consider the anti-ribbon expansion in detail.
\begin{lemma}
\label{lem:leaf-are}
Let $E$ be a leaf of the sliding tree $\mathcal{T}_D$ of a corner-crossed diagram $D$, identified with its normalization.
Then $E$ is an anti-ribbon diagram in the southeast convention.
Moreover, $\operatorname{SYT}(E)$ consists of the single standard Young tableau $T_0$ that enumerates the boxes of $E$ along its anti-ribbon path, and
\[
\operatorname{ARE}_E(A)=R_{T_0}(A)=A
\]
for every tableau $A$ on $E$, where diagrams are identified up to translation.
\end{lemma}

\begin{proof}
By Corollary~\ref{cor:leaf-antiribbon}, $E$ is an anti-ribbon diagram in the southeast convention.

Write $E=\{p_1,\ldots,p_N\}$ along the anti-ribbon path, so that $p_{k+1}$ is either the right neighbor or the lower neighbor of $p_k$.
In either case $p_k\preceq_E p_{k+1}$, so $\preceq_E$ contains the chain $p_1\prec\cdots\prec p_N$ and is a total order.
A finite totally ordered set has a unique linear extension, namely $T_0(p_k)=k$.
Finally, $k\in\operatorname{Des}(T_0)$ if and only if $p_{k+1}$ lies strictly below $p_k$, that is, if and only if the $k$-th step of the path is a downward step.
Therefore $R(\operatorname{rb}(T_0))=E$ up to translation, and $R_{T_0}(A)(p_k)=A(T_0^{-1}(k))=A(p_k)$, so $R_{T_0}(A)=A$.
\end{proof}

\begin{definition}
\label{def:slide-formal}
Let $D$ be a corner-crossed diagram with $N$ boxes and fix, at every node of the sliding tree $\mathcal{T}_D$, a consecutive admissible pair as in Section~\ref{sec:sliding}.
Each slide acts on tableaux by transporting entries along its box bijection: $(\sigma A)(\phi(\beta))=A(\beta)$ for $\beta\in D$.
For a root-to-leaf path $P$ in $\mathcal{T}_D$ ending at the leaf $D_P$, let $\Phi_P\colon D\to D_P$ denote the composite of the box bijections along $P$, and define the transported tableau $A_P$ on $D_P$ by $A_P(\Phi_P(u))=A(u)$ for $u\in D$.
The terminal sliding expansion of a tableau $A$ on $D$ is
\[
\operatorname{Slide}_D(A)
=
\sum_{P} A_P,
\]
the sum being over all root-to-leaf paths of $\mathcal{T}_D$, in the free $\mathbb{Q}$-vector space generated by anti-ribbon tableaux, where translated copies are identified.

Writing the boxes of $D_P$ in anti-ribbon path order as $p_1,\ldots,p_N$, the terminal word of $P$ is
\[
\omega(P)=s\varepsilon_1\cdots\varepsilon_{N-1},
\qquad
\varepsilon_k
=
\begin{cases}
r, & \text{if $p_{k+1}$ is immediately to the right of $p_k$},\\
b, & \text{if $p_{k+1}$ is immediately below $p_k$},
\end{cases}
\]
so that $D_P=R(\omega(P))$ up to translation.
\end{definition}

\begin{conjecture}
\label{conj:slide-are}
Let $D$ be a corner-crossed diagram without implicit boxes, so that $\widetilde{D}=D$.
Then the consecutive admissible pairs at the nodes of $\mathcal{T}_D$ can be chosen so that the terminal sliding expansion coincides with the anti-ribbon expansion:
\[
\operatorname{Slide}_D(A)
=\operatorname{ARE}_D(A)
\]
for every tableau $A$ on $D$.
In particular, after the all-one specialization, the coefficient of the anti-ribbon word $w$ in the sliding expansion is then
\[
c_D(w)
=\#\{\,T\in \operatorname{SYT}(D)\mid \operatorname{rb}(T)=w\,\}.
\]
\end{conjecture}

\begin{remark}
The restriction on implicit boxes in Conjecture~\ref{conj:slide-are} cannot simply be omitted.
Indeed, let
\[
D=\{(1,2),(2,1),(2,3)\}.
\]
Its northwest corners are $(1,2)$ and $(2,1)$, so $D$ is corner-crossed, and $(2,2)$ is an implicit box.
The unique consecutive admissible pair is $((1,2),(2,1))$.
Consider the standard tableau
\[
T=\begin{ytableau}
\none & 2 \\
1 & \none & 3
\end{ytableau}.
\]
Since $T(1,2)>T(2,1)$, its path takes the vertical slide, which produces the one-row leaf $\{(2,1),(2,2),(2,3)\}$.
Thus
\[
\omega(P_T)=\mathtt{srr}
\ne\mathtt{srb}=\operatorname{rb}(T).
\]
Both children of the root are leaves, so there are no alternative choices of pairs.
By Proposition~\ref{prop:slide-word-reduction}, the two formal expansions therefore do not agree for every tableau on this $D$.
In particular, every $M\in\mathrm{SSYT}(D)$ satisfies $M(1,2)<M(2,3)$, whereas the one-row leaf allows these two entries to be equal.
\end{remark}

\begin{proposition}
\label{prop:slide-word-reduction}
Let $D$ be a corner-crossed diagram with $N$ boxes.
For $T\in\operatorname{SYT}(D)$, let $P_T$ denote the root-to-leaf path associated with $T$ by Theorem~\ref{thm:onetoone}.
Then the following are equivalent.
\begin{enumerate}
\item $\operatorname{Slide}_D(A)=\operatorname{ARE}_D(A)$ for every tableau $A$ on $D$.
\item $\omega(P_T)=\operatorname{rb}(T)$ for every $T\in\operatorname{SYT}(D)$.
\end{enumerate}
In particular, if $D$ has no implicit boxes, Conjecture~\ref{conj:slide-are} holds for $D$ if and only if the pairs can be chosen so that the identity in (2) holds.
\end{proposition}

\begin{proof}
Fix $T\in\operatorname{SYT}(D)$ and write the boxes of $D_{P_T}$ in anti-ribbon path order as $p_1,\ldots,p_N$.
The anti-ribbon-reading property in Theorem~\ref{thm:onetoone} states that transporting $T$ along $P_T$ yields the entries $1,\ldots,N$ in this order; equivalently,
\[
\Phi_{P_T}\bigl(T^{-1}(k)\bigr)=p_k
\qquad(1\leq k\leq N).
\]
Hence, for every tableau $A$ on $D$,
\[
A_{P_T}(p_k)=A\bigl(T^{-1}(k)\bigr)
\qquad(1\leq k\leq N),
\]
so the summand $A_{P_T}$ of $\operatorname{Slide}_D(A)$ and the summand $R_T(A)$ of $\operatorname{ARE}_D(A)$ carry the same entry sequence along their anti-ribbon paths and differ at most in their shapes, which are $R(\omega(P_T))$ and $R(\operatorname{rb}(T))$ respectively.

Assume (2).
Then $A_{P_T}=R_T(A)$ for every $T$, and since $T\mapsto P_T$ is a bijection onto the set of root-to-leaf paths,
\[
\operatorname{Slide}_D(A)
=\sum_{T\in\operatorname{SYT}(D)}A_{P_T}
=\sum_{T\in\operatorname{SYT}(D)}R_T(A)
=\operatorname{ARE}_D(A).
\]

Conversely, assume (1) and apply it to the universal tableau $X(u)=x_u$ with pairwise distinct formal symbols $x_u$.
The entries of $X_{P_T}$ in path order are $x_{T^{-1}(1)},\ldots,x_{T^{-1}(N)}$, and the entries of $R_{T'}(X)$ in path order are $x_{T'^{-1}(1)},\ldots,x_{T'^{-1}(N)}$; since the symbols are distinct, these sequences agree if and only if $T'=T$.
Thus each anti-ribbon tableau occurs at most once on either side of $\operatorname{Slide}_D(X)=\operatorname{ARE}_D(X)$, and matching the terms with equal entry sequences yields $X_{P_T}=R_T(X)$, hence $R(\omega(P_T))=R(\operatorname{rb}(T))$.
Since $w\mapsto R(w)$ is injective up to translation, $\omega(P_T)=\operatorname{rb}(T)$.
\end{proof}

\section{Schur functions and Schur multiple zeta functions}
\label{sec:applications}
\subsection{Schur function}
A well-known definition of Schur polynomials is given in terms of Young tableaux.
Because the number of variables must be finite, we first restrict the entries.
For a diagram $D$ and an integer $m\ge1$, put
\[
\mathrm{SSYT}_m(D)
=\{\,T=(t_{ij})\in \mathrm{SSYT}(D)
\mid
1\leq t_{ij}\leq m \text{ for every } (i,j)\in D\,\},
\]
which is a finite set, and abbreviate
\[
\mathrm{SSYT}_m(\lambda)=\mathrm{SSYT}_m\bigl(D(\lambda)\bigr),
\quad
\mathrm{SSYT}_m(\lambda/\mu)=\mathrm{SSYT}_m\bigl(D(\lambda/\mu)\bigr),
\quad
\mathrm{SSYT}_{\mathrm{ar},m}(\alpha)=\mathrm{SSYT}_m\bigl(D_{\mathrm{ar}}(\alpha)\bigr).
\]
For a tableau $T=(t_{ij})$ on $D$ we write
\[
x^T=\prod_{(i,j)\in D} x_{t_{ij}}.
\]

\begin{definition}
\label{def:s_D}
For a diagram $D$ and variables $x_1,\ldots,x_m$, put
\[
s_D(x_1,\ldots,x_m)=\sum_{T\in \mathrm{SSYT}_m(D)} x^T .
\]
We write $s_\lambda=s_{D(\lambda)}$, $s_{\lambda/\mu}=s_{D(\lambda/\mu)}$ and $s^{\mathrm{ar}}_{\alpha}=s_{D_{\mathrm{ar}}(\alpha)}$; the first is the Schur polynomial of shape $\lambda$, the second the skew Schur polynomial of shape $\lambda/\mu$, and the third their anti-ribbon analogue attached to a composition $\alpha$.
\end{definition}

We now apply Theorem~\ref{thm:standardization}.
For a corner-crossed diagram $D$ and $T\in\mathrm{SYT}(D)$ we keep the notation $\alpha_T=\alpha(\operatorname{rb}(T))$ and $\mathcal{S}_T$ of Theorem~\ref{thm:standardization}.
By Theorem~\ref{thm:onetoone}, the map $T\mapsto P_T$ is a bijection from $\mathrm{SYT}(D)$ onto the set of root-to-leaf paths of any fixed sliding tree $\mathcal{T}_D$, and the slides along $P_T$ carry the box $T^{-1}(k)$ to the $k$-th box of the leaf in path order; hence $R_T(M)$ has the same entry sequence as the filling of the leaf obtained by transporting $M$ along $P_T$.
To obtain the semi-standard decomposition, this sequence is placed on the descent shape $D_{\mathrm{ar}}(\alpha_T)$, which coincides with the shape of the leaf exactly when $\omega(P_T)=\operatorname{rb}(T)$.

\begin{theorem}
\label{thm:schur-expansion}
Let $D$ be a corner-crossed diagram and let $m\ge1$. Then
\[
s_D(x_1,\ldots,x_m)
=\sum_{T\in\mathrm{SYT}(D)}
s^{\mathrm{ar}}_{\alpha_T}(x_1,\ldots,x_m).
\]
The number of summands, counted with multiplicity before collecting equal terms, is $|\mathrm{SYT}(D)|$.
\end{theorem}

\begin{proof}
By Theorem~\ref{thm:standardization}(1), $\mathrm{SSYT}_m(D)$ is the disjoint union of the sets $\mathcal{S}_T\cap\mathrm{SSYT}_m(D)$ with $T\in\mathrm{SYT}(D)$.
By Theorem~\ref{thm:standardization}(2), the map $M\mapsto R_T(M)$ is a bijection from each such set onto $\mathrm{SSYT}_{\mathrm{ar},m}(\alpha_T)$, and it preserves the monomial $x^M$.
Summing $x^M$ over this partition proves the identity.
\end{proof}

We illustrate Theorem~\ref{thm:schur-expansion} on Example~\ref{ex1}. Let $m=3$. The left-hand side is the Schur polynomial of shape $\lambda=(2,1)$, 
\[
s_{\lambda}(x_1, x_2, x_3)=x_1^2x_2+x_1^2x_3+x_1x_2^2+2x_1x_2x_3+x_1x_3^2+x_2^2x_3+x_2x_3^2,
\]
and the two terms on the right-hand side are the anti-ribbon analogues attached to the compositions $\alpha_1=(2,1)$ and $\alpha_2=(1,2)$, 
\[
s^{\mathrm{ar}}_{(2,1)}(x_1,x_2,x_3) = x_1^2x_2+x_1^2x_3+x_1x_2x_3+x_2^2x_3,
\] 
and
\[
s^{\mathrm{ar}}_{(1,2)}(x_1,x_2,x_3) = x_1x_2^2+x_1x_2x_3+x_1x_3^2+x_2x_3^2.
\]

We note that the polynomials of anti-ribbon type are in general not symmetric: for instance $s^{\mathrm{ar}}_{(2,1)}(x_1,x_2,x_3)$ contains the monomial $x_1^2x_2$ but not $x_1x_2^2$.
The one-row and one-column cases are exceptions, since $s^{\mathrm{ar}}_{(n)}=h_n$ and $s^{\mathrm{ar}}_{(1^n)}=e_n$.
Thus the individual terms of Theorem~\ref{thm:schur-expansion} need not be symmetric, and the decomposition is in general finer than a decomposition into Schur polynomials.
In fact, by Lemma~\ref{lem:ar-ssyt}, $s^{\mathrm{ar}}_{\alpha}$ is the fundamental quasisymmetric polynomial $F_{S(\alpha)}$ in $m$ variables, and Theorem~\ref{thm:schur-expansion} is a tableau-level form of Gessel's expansion $s_{\lambda/\mu}=\sum_{T\in\mathrm{SYT}(\lambda/\mu)}F_{\operatorname{Des}(T)}$~\cite{Gessel1984}.

\subsection{Schur multiple zeta functions and those of anti-ribbon type}
In \cite{NakasujiPhuksuwanYamasaki2018}, Nakasuji, Phuksuwan, and Yamasaki defined the Schur multiple zeta function as a zeta function which incorporates the structure of the Schur function and interpolates between the multiple zeta and the multiple zeta-star functions of the Euler--Zagier type.
\begin{definition}[Schur multiple zeta function]
 Let $\lambda$ and $\mu$ be partitions with $\mu\subset \lambda$. 
 For a set ${ \pmb s}=(s_{ij})\in T(\lambda/\mu,\mathbb{C})$ of variables, the Schur multiple zeta function associated with the shape $\lambda/\mu$ is defined by
 \[
  \zeta_{\lambda/\mu}(\pmb {s})
  \coloneqq
  \sum_{M\in \mathrm{SSYT}(\lambda/\mu)} \prod_{(i,j)\in D(\lambda/\mu)} 
  \frac{1}{m_{ij}^{s_{ij}}},
 \]
 where the sum runs over all semi-standard Young tableaux $M=(m_{ij})$ of shape $\lambda/\mu$.
\end{definition}
A sufficient condition for absolute convergence was established in \cite{NakasujiPhuksuwanYamasaki2018}; see also \cite{NakasujiTakeda2025}, where the same argument is carried out for the Schur $P$- and $Q$-analogues.
The series converges absolutely on
\[
 W_{\lambda/\mu}
=
\left\{{\pmb s}=(s_{ij})\in T(\lambda/\mu,\mathbb{C}) \mid
\begin{array}{l}
 \text{${\rm Re}(s_{ij})\ge 1$ for all $(i,j)\in D(\lambda/\mu) \setminus C(\lambda/\mu)$} \\[3pt]
 \text{${\rm Re}(s_{ij})>1$ for all $(i,j)\in C(\lambda/\mu)$}
\end{array}
\right\}.
\]
Schur multiple zeta functions have also been studied for other families of diagrams: the Schur $P$- and $Q$-multiple zeta functions of \cite{NakasujiTakeda2025} are indexed by shifted diagrams, and anti-hook shapes have been studied in \cite{MatsumotoNakasuji2021}.
The anti-ribbon shapes considered here are obtained by reflecting ribbon shapes in a vertical axis.
For each standard tableau $T$ of a skew shape, its descent set $\operatorname{Des}(T)$ determines an anti-ribbon shape.
The expansion in Theorem~\ref{thm:smzv-expansion} has a term of this shape indexed by $T$.
This motivates the following definition of Schur multiple zeta functions of anti-ribbon type.
Let $\alpha=(\alpha_1,\ldots,\alpha_r)$ be a composition.
For $\pmb{s}=(s_{ij})\in T_{\mathrm{ar}}(\alpha,\mathbb C)$, we define the Schur multiple zeta function of anti-ribbon type by
\[
\zeta^{\mathrm{ar}}_{\alpha}(\pmb{s})
=
\sum_{M=(m_{ij})\in \mathrm{SSYT}_{\mathrm{ar}}(\alpha)}
\prod_{(i,j)\in D_{\mathrm{ar}}(\alpha)}
\frac{1}{m_{ij}^{s_{ij}}}.
\]
\begin{theorem} \label{thm:ar-convergence}
Let $\alpha=(\alpha_1,\ldots,\alpha_r)$ be a composition and let
$\pmb{s}=(s_{ij})_{(i,j)\in D_{\mathrm{ar}}(\alpha)}$ be a family of complex numbers.
Let $R=\alpha_1+\cdots+\alpha_r-r+1$. 
If
\[
\mathrm{Re}(s_{ij})\ge1
\qquad ((i,j)\in D_{\mathrm{ar}}(\alpha)) \text{ and } \mathrm{Re}(s_{rR})>1,
\]
then the series
$
\zeta^{\mathrm{ar}}_{\alpha}(\pmb{s})
$
converges absolutely.
\end{theorem}

\begin{proof}
We focus on the entry at position $(r,R)$ of the tableau. Then, we have
\begin{align}
\label{for:proof_for_converge1}
\left|\sum_{\substack{M=(m_{ij})\in \mathrm{SSYT}_{\mathrm{ar}}(\alpha)\\ m_{ij}\le N}}
\prod_{(i,j)\in D_{\mathrm{ar}}(\alpha)}
m_{ij}^{-s_{ij}}\right|
\le
\sum_{N_1=1}^N
\left(
\sum_{\substack{M=(m_{ij})\in \mathrm{SSYT}_{\mathrm{ar}}(\alpha)\\ m_{ij}\le N_1}}
\prod_{\substack{(i,j)\in D_{\mathrm{ar}}(\alpha)\\ (i,j)\neq (r,R)}}
m_{ij}^{-\mathrm{Re}(s_{ij})}
\right)
\frac{1}{N_1^{\mathrm{Re}(s_{rR})}}.
\end{align}
By enlarging the region of summation, we obtain
\begin{align*}
\sum_{\substack{M=(m_{ij})\in \mathrm{SSYT}_{\mathrm{ar}}(\alpha)\\ m_{ij}\le N_1}}
\prod_{\substack{(i,j)\in D_{\mathrm{ar}}(\alpha)\\ (i,j)\neq (r,R)}}
m_{ij}^{-\mathrm{Re}(s_{ij})}
&\le
\prod_{\substack{1\le i\le r\\ 1\le j\le R\\ (i,j)\neq (r,R)}}
\sum_{m=1}^{N_1}\frac{1}{m}.
\end{align*}
Because for any $\varepsilon>0$ there exists a constant $C_\varepsilon>0$ such that
\[
\sum_{m=1}^N \frac{1}{m} \le C_\varepsilon N^\varepsilon,
\]
it follows that
\[
\sum_{\substack{M=(m_{ij})\in \mathrm{SSYT}_{\mathrm{ar}}(\alpha)\\ m_{ij}\le N_1}}
\prod_{\substack{(i,j)\in D_{\mathrm{ar}}(\alpha)\\ (i,j)\neq (r,R)}}
m_{ij}^{-\mathrm{Re}(s_{ij})}
\le
C_\varepsilon^{rR}\, N_1^{\varepsilon rR}.
\]
Combining this with \eqref{for:proof_for_converge1}, we obtain
\[
\left|\sum_{\substack{M=(m_{ij})\in \mathrm{SSYT}_{\mathrm{ar}}(\alpha)\\ m_{ij}\le N}}
\prod_{(i,j)\in D_{\mathrm{ar}}(\alpha)}
m_{ij}^{-s_{ij}}\right|
\le
C_\varepsilon^{rR}
\sum_{N_1=1}^{N}
N_1^{\varepsilon rR-\mathrm{Re}(s_{rR})}.
\]

By choosing $\varepsilon>0$ sufficiently small so that $\mathrm{Re}(s_{rR})-\varepsilon rR>1$, the right-hand side converges as $N\to\infty$. This proves the theorem.
\end{proof}
\begin{theorem}
\label{thm:smzv-expansion}
Let $\mu\subset\lambda$ be partitions with $D=D(\lambda/\mu)$, and let $\pmb{s}=(s_{\beta})\in W_{\lambda/\mu}$.
For $T\in\mathrm{SYT}(D)$, regard $R_T(\pmb{s})$ as a filling of $D_{\mathrm{ar}}(\alpha_T)$ by complex numbers, so that its entry at the $k$-th box in path order is $s_{T^{-1}(k)}$.
Then every series $\zeta^{\mathrm{ar}}_{\alpha_T}(R_T(\pmb{s}))$ converges absolutely, and
\[
\zeta_{\lambda/\mu}(\pmb{s})
=
\sum_{T\in\mathrm{SYT}(D)}
\zeta^{\mathrm{ar}}_{\alpha_T}\bigl(R_T(\pmb{s})\bigr).
\]
\end{theorem}

\begin{proof}
Put $N=|D|$ and fix $T\in\mathrm{SYT}(D)$.
The last box $p_N$ of $D_{\mathrm{ar}}(\alpha_T)$ in path order is the box denoted by $(r,R)$ in Theorem~\ref{thm:ar-convergence}, and the entry of $R_T(\pmb{s})$ there is $s_{T^{-1}(N)}$.
The box $T^{-1}(N)$ is maximal for $\preceq_D$, so neither its right nor its lower neighbor belongs to $D$, that is, $T^{-1}(N)\in C(\lambda/\mu)$.
Since $\pmb{s}\in W_{\lambda/\mu}$, every entry of $R_T(\pmb{s})$ has real part at least $1$, and the entry at $p_N$ has real part $\mathrm{Re}(s_{T^{-1}(N)})>1$.
Hence $\zeta^{\mathrm{ar}}_{\alpha_T}(R_T(\pmb{s}))$ converges absolutely by Theorem~\ref{thm:ar-convergence}.

For $M\in\mathcal{S}_T$, the boxes $T^{-1}(k)$ of $D$ and $p_k$ of $D_{\mathrm{ar}}(\alpha_T)$ carry the same entry $M(T^{-1}(k))$ and the same variable $s_{T^{-1}(k)}$, whence
\[
\prod_{\beta\in D}\frac{1}{M(\beta)^{s_\beta}}
=
\prod_{k=1}^{N}\frac{1}{M(T^{-1}(k))^{s_{T^{-1}(k)}}}
=
\prod_{\gamma\in D_{\mathrm{ar}}(\alpha_T)}\frac{1}{R_T(M)(\gamma)^{R_T(\pmb{s})(\gamma)}}.
\]
By Theorem~\ref{thm:standardization}, $\mathrm{SSYT}(\lambda/\mu)$ is the disjoint union of the sets $\mathcal{S}_T$, and $M\mapsto R_T(M)$ maps $\mathcal{S}_T$ bijectively onto $\mathrm{SSYT}_{\mathrm{ar}}(\alpha_T)$.
Since the series defining $\zeta_{\lambda/\mu}(\pmb{s})$ converges absolutely, it may be rearranged along this finite disjoint decomposition, and summing over the blocks gives the assertion.
\end{proof}

Let $\mathcal{Z}$ denote the $\mathbb{Q}$-vector space spanned by all multiple zeta values.
It is known that the $\mathbb{Q}$-vector space spanned by Schur multiple zeta values coincides with $\mathcal{Z}$~\cite{NakasujiPhuksuwanYamasaki2018}.
For $\pmb{s}\in W_{\lambda/\mu}$ with positive integer entries, Theorem~\ref{thm:smzv-expansion} expresses every Schur multiple zeta value as a finite $\mathbb{Q}$-linear combination of Schur multiple zeta values of anti-ribbon type.
Conversely, by decomposing the defining summation according to the weak and strict inequalities among the running indices, every Schur multiple zeta value of anti-ribbon type can be expressed as a finite $\mathbb{Q}$-linear combination of ordinary multiple zeta values.
Consequently, the Schur multiple zeta values of anti-ribbon type span $\mathcal{Z}$ over $\mathbb{Q}$.

\section*{Acknowledgments}
This work was supported by JSPS KAKENHI Grant Numbers JP22K13897 (Murahara), JP22K03274 (Nakasuji), and JP22K13900 (Takeda).

\bibliographystyle{plain}
\bibliography{References} 

\end{document}